\documentclass{article}

\usepackage[latin1]{inputenc}
\usepackage[T1]{fontenc}
\usepackage{amsmath}
\usepackage{amssymb}
\usepackage{theorem}
\usepackage{amscd}

\newtheorem{thm}{Theorem}
\newtheorem{defn}{Definition}
\newtheorem{lemma}{Lemma}

\newtheorem{cory}{Corollary}
\newtheorem{rmk}{Remark}

\newcommand{\oct}{\mathbb{O}}
\newcommand{\R}{\mathbb{R}}
\newcommand{\C}{\mathbb{C}}
\newcommand{\h}{\mathbb{H}}
\newcommand{\g}{\mathfrak{g}}
\newcommand{\gt}{\tilde{\mathfrak{g}}}
\newcommand{\lm}{\lambda}
\newcommand{\Lm}{\Lambda}
\newcommand{\iif}{if and only if }

\newcommand{\im}{\text{\normalfont Im\,}}
\newcommand{\eps}{\varepsilon}
\newcommand{\G}{\mathcal{G}}
\newcommand{\dl}[2]{\frac{\partial{#1}}{\partial{#2}}}
\newcommand{\mk}{\mathfrak{m}}
\newcommand{\hk}{\mathfrak{h}}
\newcommand{\lie}{\mathrm{Lie}\,}
\newcommand{\Is}{\mathrm{Is}}
\newcommand{\Ad}{\mathrm{Ad}}
\newcommand{\ad}{\mathrm{ad}}
\newcommand{\so}{\mathfrak{so}}
\newcommand{\mrm}[1]{\mathrm{#1}}

\newcommand{\tm}{\tau_{|\mathfrak{m}}}
\newcommand{\Int}{\mathrm{Int}}
\newcommand{\iso}{\mathrm{Is}_{p_0}}
\newcommand{\tl}[1]{\tilde{#1}}
\newcommand{\aut}{\mathrm{Aut}}
\newcommand{\der}{\mathrm{Der}}
\newcommand{\unt}{\underline{\tau}}
\newcommand{\taum}{\tau_{\mathfrak{m}}}
\newcommand{\tauh}{\tau_{\mathfrak{h}}}
\newcommand{\tmh}{\tau_{\mid\mathfrak{h}}}
\newcommand{\und}[1]{\underline{#1}}
\newcommand{\mak}[1]{\mathfrak{#1}}
\newcommand{\Adm}{\mathrm{Ad}_{\mk}}
\newcommand{\adm}{\mathrm{ad}_{\mk}}
\newcommand{\Id}{\mathrm{Id}}
\newcommand{\Z}{\mathbb{Z}}
\newcommand{\diag}{\mathrm{diag}}
\newcommand{\lag}{\mathrm{Lag}}
\newcommand{\lis}{\mathfrak{Is}}
\newcommand{\liso}{\mathfrak{Is}_{p_0}}

\author{Idrisse Khemar}
\title{Geometric Interpretation of Second Elliptic Integrable Systems}
\date{}

\begin{document}
\maketitle

\null\hfill  \textbf{Abstract.}  \hfill\null\\
In this paper we give a geometrical interpretation of all the second elliptic integrable
systems associated to 4-symmetric spaces. We first show that a 4-symmetric space $G/G_0$ can
be embedded into the twistor space of the corresponding symmetric space $G/H$. Then we
prove that the second elliptic system is equivalent to the vertical harmonicity of an
admissible twistor lift $J$ taking values in $G/G_0 \hookrightarrow \Sigma(G/H)$. We begin the paper with
an example: $G/H=\R^4$. We also study  the structure of 4-symmetric bundles over Riemannian symmetric spaces.\\[1mm]
\section*{Introduction}
The first example of second elliptic integrable system associated to a 4-symmetric space was given
in \cite{HR1}: the authors showed that the Hamiltonian stationary Lagrangian surfaces in $\C^2$ are solutions of one such
integrable system. Later they generalized their result to complex two-dimensional Hermitian symmetric spaces, \cite{HR3}.
In \cite{ki1}, we presented a new class of geometric problems for surfaces in the Euclidean space of dimension 8 by  identifying $\R^8$ with the set of octonions $\oct$, and we proved that these problems are solutions of a second elliptic integrable system. Using the left multiplication in $\mathbb{O}$ by the vectors of the canonical basis of $\im\oct$  we defined
 a family $\{\omega_{i},\,1\leq i\leq 7\}$ of canonical symplectic forms in $\oct$. This allowed us to define
 the notion of \emph{$\omega_{I}$-isotropic surfaces}, for $I\varsubsetneqq\{1,...,7\}$. 
  Using the cross-product in $\mathbb{O}$ we defined a map $\rho\colon
Gr_2(\mathbb{O})\to S^6$ from the Grassmannian of planes in $\mathbb{O}$ to
$S^6$. This allowed us to associate to each surface $\Sigma$ in $\mathbb{O}$
a function $\rho_{\Sigma}\colon \Sigma\to S^6$. In the case of $\omega_{I}$-isotropic surfaces, $\rho_{\Sigma}$
takes values in a subsphere $S^{I}=S(\oplus_{i \notin I,i>0}\R e_i)\simeq S^{6-|I|}$.
 We showed that
the  surfaces in $\mathbb{O}$ such that $\rho_{\Sigma}$ is
harmonic (\emph{$\rho$-harmonic surfaces}) are solutions of a completely integrable system $\mathcal{S}$. 
More generally we showed that the $\omega_{I}$-isotropic $\rho$-harmonic surfaces are solutions of a 
completely integrable system $\mathcal{S}_{I}$. Hence we built a family $(\mathcal{S}_I)$ indexed by $I$,
  of  set of surfaces solutions of an integrable system, all included in
 $\mathcal{S}=\mathcal{S}_\emptyset$, such that $I\subset J$ implies $\mathcal{S}_J
\subset\mathcal{S}_I$. Each $\mathcal{S}_I$ is a second elliptic integrable system (in the sense of C.L. Terng).
 This means that the equations of this system are equivalent to the zero curvature equation :
$$
d\alpha_{\lm} + \frac{1}{2}[\alpha_{\lm}\wedge\alpha_{\lm}]=0,
$$
for all $\lambda\in \C^*$, and where $\alpha_{\lm}=\lm^{-2}\alpha_2' + {\lm}^{-1}\alpha_{-1} + \alpha_0 +
\lm\alpha_{1} + \lm^2\alpha_2''$.\\
By restriction to the quaternions $ \mathbb{H}\subset\mathbb{O}$ of our theory we obtain   
a new class of surfaces: the $\omega_{I}$-isotropic $\rho$-harmonic surfaces in $\h$. Then $\rho(Gr_2(\h))=S^2$
and $|I|=0,1$ or $2$. For $|I|=1$ we obtain the Hamiltonian Stationary Lagrangian surfaces in $\R^4$ 
and for $|I|=2$, the special  Lagrangian surfaces. By restriction to $\im\h$, 
we obtain the CMC surfaces of $\mathbb{R}^3$.\\[0.1cm]
Besides, in \cite{ki2}, we found a supersymmetric interpretation of all the second elliptic integrable systems associated to  4-symmetric spaces in terms of super harmonic maps into  symmetric spaces. This led us to conjecture that this system has a geometric interpretation in terms of surfaces with values in a symmetric space, such that a certain associated map is harmonic as this is the case for Hamiltonian stationary
Lagrangian surfaces in Hermitian symmetric spaces 
 or for $\rho$\,-harmonic surfaces of $\oct$.\\
In this paper we give the answer to this conjecture. More precisely, 
we give a geometric interpretation -- in terms of  vertical harmonic twistor lifts  -- of all the
second elliptic integrable systems associated to  4-symmetric spaces. 
Indeed given a 4-symmetric space $G/G_0$, and its  order four automorphism
$\tau\colon G\to G$, then the involution $\sigma=\tau^2$ gives rise to the symmetric space
$G/H$, with $H=G^{\sigma}$. Then we prove that the 
second elliptic integrable system associated to the 4-symmetric space $G/G_0$ is exactly the
equation of vertical harmonicity for an admissible twistor lift in $G/H$.
More  precisely, given a 4-symmetric space $G/G_0$, and its associated symmetric space $G/H$, then $G/G_0$ is a subbundle of the twistor space $\Sigma(G/H)$. 
We  prove that the second elliptic integrable systems associated to $G/G_0$, is the system of equations for maps $J\colon\C\to G/G_0\subset\Sigma(G/H)$ such that
$J$ is compatible with the Gauss map of $X\colon\C\to G/H$, the projection of $J$ into $G/H$, i.e. $X$ is $J$-holomorphic
(admissible twistor lift), and such that $J$ is vertically harmonic. We prove also that an admissible twistor lift $J\colon\C\to G/G_0$ is harmonic \iif it is vertically harmonic and $X\colon \C\to G/H$ is harmonic.\\
 We begin the paper with an example: $\R^4$. This case was just mentioned briefly at the end of \cite{ki1} as a restriction of the difficult problem in $\oct$. In this paper we study this problem independently and in detail. However,  we  also present a formulation of this problem in terms of twistor lifts which seems to be the appropriate formulation. Besides, in dimension 4 we have unicity of the twistor lift (in $\Sigma^+(G/H)$ and $\Sigma^-(G/H)$ respectively) so we are in this case in the presence of a theory of surfaces (and not, as in the general case, a theory of twistor lift). Hence we can speak about $\rho$-harmonic surfaces in this dimension (which are exactly the solutions of the second elliptic integrable system).
In our work we are led to prove some theorems on the structure of 4-symmetric bundles. Indeed we want to answer  the following questions. Given a Riemannian symmetric space, do there exist  4-symmetric bundles over it? In other words, does its twistor bundle contain 4-symmetric subbundles, and if yes, how can we characterize these 4-symmetric components? are they isomorphic? How are they distributed in the twistor space ? Do they form a partition of the twistor space? 
The 4-symmetric spaces have been classified (at least in the compact case, see \cite{jimenez,wolf}). However, our point of view is different: we want to  keep an intrinsic point of view as long as possible, therefore we deal with the Riemannian symmetric space and a (locally) 4-symmetric bundle defined over it, and we try to forget as much as possible the order four automorphism of the Lie algebra. Our aim is to give a formulation of our problem which is as general and intrinsic as possible. For example, our definition of vertical harmonicity holds for any Riemannian manifold. Moreover we prove the following characterization: to define a (locally) 4-symmetric bundle over $M$ is equivalent to give ourself $J_0\in\Sigma(T_{p_0}M)$, an (orthogonal) almost complex structure in $T_{p_0}M$, which leaves invariant the curvature. We obtain the following picture: the submanifold of the twistor bundle leaving invariant the curvature is the disjoint union of all the maximal (locally) 4-symmetric subbundle, which are orbits (under the action of some subgroups of $\Is(M)$). Each isomorphism class of orbits defines a different second elliptic integrable system.\\
Our paper is organized as follows. In Section~\ref{r4} we deal with the $\rho$-harmonic surfaces in $\R^4$. Section~\ref{heart} contains our main result: the interpretation of the second elliptic integrable systems associated to a 4-symmetric space in terms of  vertical harmonicity of an admissible twistor lift. Then Sections~\ref{structure} and \ref{splitting} are devoted to the study of the structure of 4-symmetric bundles over symmetric spaces. The last Section presents some examples of 4-symmetric bundles.\\
%
%
%\footnote{In all the paper, $\Sigma(M)$ denote the twistor space of the manifold $M$.}
%  
\section{$\rho$-harmonic surfaces in $\h$}\label{r4}
\subsection{Cross product, complex structure and Grassmannian of
planes in $\h$}

We consider the space $\R^4=\h$ with its canonical basis $(1,i,j,k)$
(which we denote also by $(e_i)_{0\leq i\leq 3}$). Let $P=q\wedge q'$
be an oriented plane of $\h$ (itself oriented by its canonical basis)
then there exists an unique positive  complex structure\footnote{In all the paper, for any oriented Euclidean space $E$, $\Sigma(E)=\{J\in SO(E)|J^2=-\Id\}$, and $\Sigma(M)$ denotes the twistor bundle of the Riemannian manifold $M$.}
$I_P\in\Sigma^{+}(P)$ on the plane $P$. It is defined by $I_P(q)=
q'$, $I_P(q')=-q$ if $(q,q')$ is orthogonal. Next, we can  extend it
in an unique way to a positive (resp. negative) complex structure in
$\h=P\oplus P^{\bot}$, $J_{P}^+$ (resp. $J_{P}^-$) given
by
\begin{eqnarray}
J_{P}^+ & = & I_P\oplus I_{P^{\bot}} \nonumber \\
J_{P}^- & = & I_P\oplus -I_{P^{\bot}}
\end{eqnarray}
($P^{\bot}$ is oriented  so that $=P\oplus P^{\bot}$ is positively oriented).
Hence we obtain a surjective map:
\begin{equation}\label{j+}
  \begin{array}{crcl}
   J^+\colon & Gr_2(\h) & \to & \Sigma^+(\h)\\
     & q\wedge q' & \mapsto  & J_{q\wedge q'}^+  \end{array}
\end{equation}
$Gr_2(\h)$ being the Grassmannian of oriented planes in $\h$, and in the same way a  surjective map $J^-\colon Gr_2(\h)\to\Sigma^-(\h)$.\\
Besides, we have
$$J_{q\wedge q'}^+ = L_{q{\times}_L q'}=\frac{1}{2}(L_{q'}L_{\overline{q}}
 - L_qL_{\overline{q'}}),$$
 where $q{\times}_L q'=-\im(q\cdot\overline{q'})$ $=\im(q'\cdot
\overline{q})$ is the left  cross product (it is a bilinear skew map from
$\h\times\h$ to $\im\h$). Indeed, if $(q,q')$ is orthonormal then  $q{\times}_L
q'=-q\cdot\overline{q'}\in S(\im\h)$ so $L_{q{\times}_L q'}$ is a
complex structure in $\h$ and it is  positive (because $\{L_u,u\in
S^2\}$ is connected and $L_i\in\Sigma^+(\h)$ because $(1,L_i(1),j,L_i(j))$
is positively oriented). Moreover if $(q,q')$ is orthonormal then $L_{q{\times}_L q'}
(q)=(q'\overline q)q=q'$. Hence $L_{q{\times}_L q'}=J_{q\wedge q'}^+$.
Thus we obtain a diffeomorphism:
\begin{equation}\label{j(1)}
  \begin{array}{rcl}
  \Sigma^+(\h)  & \xrightarrow{\sim}& S^2\\
    J  & \longmapsto & J(1)
  \end{array}.
\end{equation}
Under this identification, the map (\ref{j+}) becomes
$$
 \begin{array}{crcl}
\rho_+ \colon & Gr_2(\h) & \to & S^2 \\
     & q\wedge q' & \mapsto  & q\times_{L}q' \ . \end{array}
$$
We can do the same for $\Sigma^-(\h)$. We obtain that
$J_{q\wedge q'}^- = R_{\overline{q{\times}_R q'}}=-R_{q{\times}_R q'}=
\frac{1}{2}(R_{q'}R_{\overline{q}} - R_qR_{\overline{q'}})$, where
$q{\times}_R q'=-\im(\overline q\cdot q')$ $=\im(\overline{q'}\cdot q)$ is the
right  cross product (it is a bilinear skew map from
$\h\times\h$ to $\im\h$). Then we have the same identification
between $\Sigma^-(\h)$ and $S^2$, as in (\ref{j(1)}). Under this
identification $J^-$  becomes
$$
 \begin{array}{crcl}
\rho_- \colon & Gr_2(\h) & \to & S^2 \\
     & q\wedge q' & \mapsto  & q\times_{R}q' \ . \end{array}
$$
\subsection{Action of $SO(4)$}
Recall the following 2-sheeted covering of $SO(4)$:
$$
  \begin{array}{crcl}
\chi\colon & S^3\times S^3 &\to & SO(4)\\
 & (a,b) &\mapsto & L_aR_{\overline b}\end{array}
$$
and set $Spin(3)_+ =\{L_a,a\in S^3\}$, $Spin(3)_-=\{R_{\overline b},b\in
S^3\}$, then $SO(4)=Spin(3)_+ Spin(3)_- =Spin(3)_- Spin(3)_+ $.
We have the two following representations of $Spin(3)_{\varepsilon}$:
$$\chi^+\colon L_a\mapsto \mathrm{int}_a =L_aR_{\overline a}\in SO(\im\h),\quad
\chi^-\colon R_{\overline b}\mapsto \mathrm{int}_b = L_bR_{\overline b}\in SO(\im\h).$$
 Then the map $\rho_{ \varepsilon}$ is $Spin(3)$-equivariant: for all
 $q,q'\in\h$, $g=L_aR_{\overline b}\in SO(4)$,
 $$
  \begin{array}{cllll}
   (gq)\times_L(gq') & = & a(q\times_L q')\overline a & = &
   \mathrm{int}_a(q\times_L q')\\
   (gq)\times_R(gq') & = & b(q\times_R q')\overline b & = &
   \mathrm{int}_b(q\times_R q').
  \end{array}$$
Hence we have $\forall g\in SO(4)$,
$$\rho_\varepsilon(g(q\wedge q'))=\chi_g^\varepsilon(\rho_\varepsilon
(q\wedge q'))$$
(where we have extended $\chi^{\varepsilon}$ to $SO(4)$ in an obvious
way: $\chi^+(L_aR_{\overline b})=\chi^+(L_a)$, $\chi^-(L_aR_{\overline b})=
\chi^-(R_{\overline b})$). Besides the map $J^\varepsilon$ is also
$Spin(3)$-equivariant, in other words the identification
(\ref{j(1)}) is $Spin(3)$-equivariant:
$$\begin{array}{lr}
\forall g\in SO(4), & \\
 &  gJ_{q\wedge q'}^+ g^{-1}= L_aR_{\overline
b}\, L_{q\times_L q'}\, R_bL_{\overline a}=L_{a(q\times_L q')a^{-1}}=J_{g(q\wedge
q')}^+  \ .\end{array}$$
The action of $Spin(3)_+=SU(\R^4,R_e)$ (resp. $Spin(3)_-=SU(\R^4,L_e)$) on
$\Sigma^-(\h)$ (resp. $\Sigma^+(\h)$) is trivial. Hence $SO(4)$ acts
on $\Sigma^\varepsilon(\h)$ only by its component
$Spin(3)_\varepsilon$ (in the same way it acts on $S_\varepsilon^2$
only by its component $Spin(3)_\varepsilon$ via
$\chi^\varepsilon $). In fact, the equality $ gJ_{q\wedge q'}^+ g^{-1}=
J_{g(q\wedge q')}^+$ results immediately from the definition of
$J_{q\wedge q'}^+$ and the fact that $g$ is a positive isometry.
This natural equality which is equivalent to what we called the
fundamental property in \cite{ki1}: $(gq)\times (gq')=\chi_g(q\times
q')$, is characteristic of  dimension 4: in this case it is
possible to associate in a natural way (which depends only on the
metric and the orientation) to each plane a complex structure, which is not possible in higher dimension.
 In dimension 8, we must choose an octonionic structure  in $\R^8$ to do that (see \cite{ki1}).

\subsection{The Grassmannian $Gr_2(\h)$ is a product of spheres}
\begin{thm}\label{product}
The map
$$  \begin{array}{crcl}
 \rho_{+}\times\rho_-\colon & Gr_2(\h) & \to & S^2\times S^2\\
 &  q\wedge q' & \mapsto & ( q\times_L q',q\times_R q') \end{array}$$
is a diffeomorphism.
\end{thm}
\textbf{Proof.} $SO(3)\times SO(3)$ acts transitively on $S^2\times
S^2$ so $SO(4)$ acts transitively on $S^2\times S^2$ via
$\chi^+\times\chi^-$, thus $\rho_+\times\rho_-$ is surjective.\\
Let $e\in S(\im\h)$, $g=L_aR_{\overline b},\, g'=L_{a'}R_{\overline{b'}}\in
SO(4)$ then we have\footnote{setting $S^1(e)=\{\cos\theta + \sin\theta\, e,\,\theta \in \R\}$,}
$$
\begin{array}{lcl}
\rho_+\times\rho_-(g(1\wedge e))=\rho_+\times\rho_-(g'(1\wedge e)) &
   \Longleftrightarrow &  (aea^{-1},-beb^{-1})=(a'e{a'}^{-1},-b'e{b'}^{-1})\\
    & \Longleftrightarrow & {a'}^{-1}a, {b'}^{-1}b\in S^1(e)\\
     & \Longrightarrow & (L_{a'}R_{\overline{b'}})^{-1}(L_{a}R_{\overline{b}})
     (1\wedge e)= 1\wedge e \\
     & \Longrightarrow  & g(1\wedge e)=g'(1\wedge e).   \end{array}$$
Hence, since $SO(4)$ acts transitively on $Gr_2(\h)$, we have proved
that $\rho_+\times\rho_-$ is injective and that
$$
\rho_+\times\rho_-(g(1\wedge e))=\rho_+\times\rho_-(g'(1\wedge e))
\Longleftrightarrow ({a'}^{-1}a, {b'}^{-1}b)\in S^1(e)\times S^1(e)
$$
(in the previous sequence of implications, the last proposition implies
the first one so all the propositions are equivalent). This completes
the proof.\hfill$\blacksquare$\\
As it is the case in \cite{ki1}, it is useful here to introduce a
function $\tilde{\rho}_\varepsilon$ on $Spin(3)_\varepsilon$
corresponding to $\rho_\varepsilon$: we define
 $\hbox{$\tilde{\rho}_\varepsilon$}_e\colon Spin(3)_\varepsilon\to
 S^2$ by
 $\hbox{$\tilde{\rho}_\varepsilon$}_e(g)=\chi_g^\varepsilon (e)$
(where $e\in S(\im\h)=S^2$), i.e. under the identification $Spin(3)_\varepsilon
=S^3$ we have $\hbox{$\tilde{\rho}_\varepsilon$}_e
(a)=\mathrm{int}_a(e)=aea^{-1}$, which is nothing but the  Hopf
fibration $S^3\to S^3/S^1(e)$. If $\rho_\varepsilon(e_1\wedge
e_2)=e$ then $\hbox{$\tilde{\rho}_\varepsilon$}_e (g)=
\rho_\varepsilon(g(e_1\wedge e_2))$. In the following, we will
forget the index $e$.\
Hence, if we take $e_1\wedge e_2$ such that $\rho_\varepsilon(e_1\wedge
e_2)=e$ for $\varepsilon=\pm 1$ (i.e. $e_1\wedge e_2=(1\wedge e)^\bot$ which
means also that $(e,e_1,e_2)$ is a direct orthonormal basis of $\im\h$)
then we have the following commutative diagram:
$$
\begin{CD}
S^3\times S^3 @> \chi >> SO(4)\\ @V {\tilde{\rho}_+\times\tilde{\rho}_-} VV
@VV {\underset{g(e_1\wedge e_2)}{\overset{g}{\downarrow}} } V \\
S^2\times S^2 @< \simeq < \rho_+\times\rho_- < Gr_2(\h)
\end{CD}
$$
Let us now consider the restriction  to $\im\h=\R^3$ of this
diagram. First the universal covering $Spin(3)\to SO(3)$ is
obtained by restriction to $\Delta_3=\{(a,a),\,a\in S^3\}\simeq S^3$
of $\chi\colon S^3\times S^3 \to SO(4)$, which gives  the covering
$(a,a)\mapsto \mathrm{int}_a$. Then supposing in addition that
$e_1,e_2\in\im\h$, the restriction  to $SO(3)$ of $SO(4)\to
Gr_2(\h)$ is only the surjective map $g\in SO(3)\mapsto g(e_1\wedge
e_2)\in Gr_3(\R^3)$. And the restriction to $Gr_2(\R^3)$ of
$\rho_+\times\rho_-$ gives the diffeomorphism $\rho\colon u\wedge v\in
Gr_2(\R^3)\to u\times v\in S^2$. Finally the restriction to
$\Delta_3$ of $\tilde{\rho}_+\times\tilde{\rho}_-$ gives  the Hopf
fibration $\tilde{\rho}\colon a\in S^3\mapsto aea^{-1}\in S^2$. So
by restriction to $\R^3$, we obtain the classical commutative
diagram:
$$
\begin{CD}
S^3 @> \chi_3 >> SO(3)\\ @V {\mrm{Hopf}} VV
@VVV \\ S^2  @< \simeq << Gr_2(\R^3)
\end{CD}
$$
\begin{rmk}\emph{
Besides if we use ${\Sigma}^{\varepsilon}(\h)$ instead of the sphere $S^2$
the Hopf fibration $\tilde{\rho}_\varepsilon$ becomes $SU(2,J_{1\wedge e}^{-
\varepsilon})\to \Sigma^{\varepsilon}(\h)=SU(2,J_{1\wedge e}^{-\varepsilon})/U(1)_\varepsilon=
SO(4)/U(2,J_{1\wedge e}^{\varepsilon})$ where
$U(1)_+=R_{S^1(e)}=\exp(\R .R_e)$, $U(1)_-= L_{S^1(e)}=\exp(\R .L_e)$.}
\end{rmk}

\subsection{The $\rho$-harmonic $\omega_{I}$-isotropic surfaces}
We recall here in the particular case of $\h=\R^4$ our result
obtained in \cite{ki1} about $\rho$-harmonic surfaces. To do that, we
need to introduce some notations and definitions. We have
$$\rho_\varepsilon(q\wedge q')=-\eps\sum_{i=1}^{3}
\omega_i^{\varepsilon}(q,q') e_i $$
where $(e_i)_{1\leq i\leq 3}=(i,j,k)$ and
$\omega_i^{\varepsilon}=\langle \cdot,J_{1\wedge e_i}^\varepsilon\cdot\rangle
$ (i.e. $\omega_i^+=\langle \cdot,L_{e_i}\cdot\rangle$, $\omega_i^-=\langle
 \cdot,R_{e_i}\cdot\rangle$). Let us set, for $I\varsubsetneqq\{1,2,3\}$,
$$
Q_I^\eps=\{P\in Gr_2(\h)|\,\omega_i^{\eps}(P)=0,\,i\in I\},
$$
then $Q_{\varnothing}=Gr_2(\h)$, $Q_{\{k\}}=\{P\in Gr_2(\h),\text{ Lagrangian for
}\omega_k^\eps\}$, and $Q_{\{k,l\}}^\eps$ is the set of special Lagrangian planes
(more precisely the $\omega_k^\eps$-Lagrangian planes  $P$  such that $\det_{\C^2}(P)=
\pm i$ under the identification: $x\in \R^4\mapsto (x_0 + i x_k, x_l
+i\eps x_{k\wedge l})\in \C^2$, with $(k,l,k\wedge l)$ cyclic permutation of
$(1,2,3)$; for example, if $(k,l)=(1,2)$, it is the identification $(z_1,z_2)\in \C^2
\mapsto z_1+z_2j\in \h$ for $\eps=1$ and $(z_1,z_2)\mapsto z_1+jz_2$ for $\eps=-1$).
We have also $\rho_{\eps}(Q_I)=S^{I}=S(\bigoplus_{i\notin I}\R e_i)=S^2,S^1,\{\pm e_k\}$
for $|I|=0,1,2$ respectively. Besides we have for $I=\{i\}\subset\{1,2,3\}$, that
$J^{+}(Q_{I})=L_{S^{I}}=S^1(\R L_{e_j}\oplus\R L_{e_k})$ is the circle of positive
complex structures which anticommute with $L_{e_i}$; and for $I=\{i,j\}\subset \{1,2,3\}
$, $J^{+}(Q_{I})=L_{S^{I}}=\{\pm L_{e_k}\}$. \\
We denote by
 $G_I^\varepsilon$ the subgroup of $Spin(3)_\varepsilon$ which conserves
 $\omega_i^\varepsilon$, for all $i\in I$; this is the subgroup of $Spin(3)_\varepsilon$ which
commutes  with $L_{e_i}$, for all $i\in I$. Then
 $G_I^\varepsilon=S^3,S^1,\{\pm 1\}$ for $|I|=0,1,2$ respectively. We
 can also consider instead of $Spin(3)_\varepsilon$ the group
 $SO(4)$ (which is equivalent to add the component  $Spin(3)_{-\varepsilon}$
 which is useless), then we have $G_{I}^\varepsilon=SO(4),
 U(2,J_{1\wedge e_i}^\varepsilon), SU(2,J_{1\wedge e_i}^\varepsilon)$ for $|I|=0,1,2$
 respectively. Let $e\in S(\bigoplus_{i\notin I}\R e_i)$.
The  inner automorphism, $\mathrm{Int} J_{1\wedge
 e}^\varepsilon$, defines on  $G_I^\varepsilon$  an involution which gives rise to the symmetric space
$S^I=G_I^\varepsilon/G_{I\cup\{k\}}^\varepsilon$ and in 
 the Lie algebra of $G_I^\varepsilon$, $\g_I^\varepsilon$, to the eigenspace decomposition of $\mathrm{Ad}J_{1\wedge e}^\varepsilon$:
$$\
\g_I^\eps=\g_0^\eps(I)\oplus\g_2^\eps(I)
$$ 
with
$\g_0^\eps(I)=\ker(\mathrm{Ad}J_{1\wedge e}^+ - \mathrm{Id})$,
$\g_2^\eps(I) =\ker(\mathrm{Ad}J_{1\wedge e}^\varepsilon
+\mathrm{Id})$.\\ Let us introduce
$\mathcal{G}_I^\eps=G_I^\eps\ltimes\R^4$ the group of affine
isometries of which the linear part is in $G_I^\eps$, and its Lie
algebra: $\g^\eps(I)=\g_I^\eps\oplus\R^4$. Consider the automorphism
of the group $\mathcal{G}_I^\eps$: $\tau_e^\eps=\mathrm{Int}(-\eps
J_{1\wedge e}^\eps,0)$ with $e \in S(\bigoplus_{i\notin I}\R e_i) $. This is an order
four automorphism which gives us an eigenspace decomposition of
$\g^\eps(I)^\C$:
$$\g^\eps(I)^\C=\bigoplus_{k\in\mathbb{Z}_4}\gt_k^\eps(I)$$ with
$\gt_{\pm 1}^\eps(I)=\g_{\pm 1}^\eps=\ker(J_{1\wedge e}^\eps\pm
i\mathrm{Id})$, $\gt_0^\eps(I)=\g_0^\eps(I)^\C$, $\gt_2^\eps(I)=\g_2^\eps(I)^\C$.
Moreover we have $[\gt_k^\eps(I),\gt_l^\eps(I)]\subset\gt_{k+l}^\eps(I)$.\\
We fix a value of $\eps=\pm 1$. Then let us define as in \cite{ki1}:

\begin{defn}
Let $L$ be an immersed surface in $\h$, then a map $\rho_{L}\colon L
\to S^2$ is associated to it, defined by $\rho_{L}(z)=\rho_{\eps}(T_z L)$
i.e. if $X\colon L\to \h$ is the immersion then $\rho_{L}=X^*\rho_\eps$.
We will say that $L$ is $\rho$-harmonic if $\rho_L$ is
harmonic\footnote{with respect to the induced metric on $S^2$.}.\\
 Let $I\varsubsetneq\{1,2,3\}$, we will say that
$L$ is $\omega_I$-isotropic if $\forall z \in L,\,T_z L\in
Q_I^\eps$ (i.e. $L$ is $\omega_i^\eps$-isotropic for all $i\in I$) which is equivalent to:
  $\rho_{L}$ takes values in $S^{I}=S(\oplus_{i\notin
 I}\R e_i)\subset S^2$. Hence for $|I|=1$, the $\rho$-harmonic $\omega_I^\eps$-isotropic
surfaces are the \textbf{Hamiltonian stationary Lagrangian surfaces} in $\C^2$, and for
$|I|=2$, these are the \textbf{special Lagrangian surfaces} in $\C^2$ (see
above for the identification $\R^4\simeq\C^2$).\\
 If it could be an ambiguity as concerned the value of $\eps=\pm
 1$, we will use the qualificatifs "left" and "right" respectively to
 design these two values. A lifted conformal left (resp. right)
 $\omega_I$-isotropic immersion - LC$\omega_I$ -
 (if $I=\varnothing$  we will say  a lifted conformal immersion or simply a lift)
 is a map $U=(F,X)\colon L\to\G_I^\eps$ such that  $X$ is a conformal
 $\omega_I^\eps$-isotropic immersion and $\tilde{\rho}_e\circ
 F=\rho_{L}$.
\end{defn}

We have obtained the following result in \cite{ki1}:

\begin{thm}
Let $\Omega$ be a simply connected open domain in $\C$, and $\alpha$ an 1-form on
$\Omega$ with values in $\g(I)$, then
\begin{description}
  \item[$\bullet$] $\alpha$ is the Maurer-Cartan form of a   LC$\omega_I$ \iif
$$d\alpha + \alpha\wedge\alpha =0, \quad \alpha_{-1}''=0\quad
\text{and}\quad \alpha_{-1}'\text{ does not vanish}$$
  \item[$\bullet$]furthermore, $\alpha$ corresponds to a
  $\rho$-harmonic $\omega_I$-isotropic conformal immersion \iif the
  extended Maurer-Cartan form $\alpha_{\lambda}=\lambda^{-2}\alpha_2' +
  \lambda^{-1}\alpha_{-1} +\alpha_0 + \lambda\alpha_1 +
  \lambda^2\alpha_2''$ satisfies
  $$d\alpha_{\lambda} + \alpha_{\lambda}\wedge\alpha_{\lambda}=0,\
\forall \lambda\in \C^*.$$
  \end{description}
\end{thm}
Let us recall the proof given in \cite{ki1}.\\[1.5mm]
\textbf{Proof.} To fix ideas, we take $\eps=1$.  $\alpha$ is a
Maurer-Cartan form \iif it satisfies the Maurer-Cartan equation. In
this case, it can be integrated by $U=(F,X)\colon\Omega\to\G_I$,
$\alpha=U^{-1}.dU,\,U(z_0)=1$. Hence $\alpha=U^{-1}.dU=(F^{-1}.dF,F^{-1}.dX)$.
Moreover, $F^{-1}.dX=\alpha_{-1} +\alpha_1$ is real and $\g_{\pm
1}=\{V\pm iL_e V,\,V\in\h\}$ so $\alpha_{-1}=\overline{\alpha_1}$.
Hence $\alpha_{-1}''=0 \Longleftrightarrow \alpha_{-1}''=\overline{\alpha_1'}=0\Longleftrightarrow
\alpha_{-1} =(F^{-1}\dl{X}{z})dz\Longleftrightarrow F^{-1} \dl{X}{y}=L_e(F^{-1}
 \dl{X}{x})\Longleftrightarrow F^{-1}dX=h(q_0du + q_0' dv)$ with $h\in
C^{\infty}(\Omega,\R)$, $q_0,q_0'\in C^{\infty}(\Omega,S^3)$,
$\langle q_0,q_0'\rangle=0$ and $\rho(q_0,q_0')=e$. Thus we have
($\alpha_{-1}''=0$ and $\alpha_{-1}'\neq 0)\Longleftrightarrow dX =e^f
(qdu + q' dv)$ with $f\in C^{\infty}(\Omega,\R)$, $(q,q')$ orthonormal and
$\rho(q,q')=\tilde{\rho}_e(F)$ i.e. $\rho_X=\tilde{\rho}_e(F)$. This
proves the first point.\\
Hence we have   the decomposition $\alpha=\alpha_2 +\alpha_{-1} +\alpha_0 +\alpha_1=
\alpha_2' +\alpha_{-1}' +\alpha_0 +\alpha_1'' +\alpha_2''$.
Furthermore, using the commutation relations $[\gt_k(I),\gt_l(I)]\subset
\gt_{k+l}(I)$, $[\g_{\pm 1},\g_{\pm 1}]=\{0\}$,  we obtain
\begin{eqnarray*}
d\alpha_{\lambda} + \alpha_{\lambda}\wedge\alpha_{\lambda} & = &
  \lambda^{-2}(d\alpha_2' + [\alpha_0\wedge\alpha_2')])\\
 &  & +\lambda^{-1}(d\alpha_{-1}' + [\alpha_{-1}'\wedge\alpha_0] +
[\alpha_1''\wedge\alpha_2']) \\
 &   &  +(d\alpha_0 + \frac{1}{2}[\alpha_0\wedge\alpha_0] + \frac{1}{2}[\alpha_2'\wedge
\alpha_2''])\\
 &  & + \lambda(d\alpha_{1}'' + [\alpha_{1}''\wedge\alpha_0] +
 [\alpha_{-1}'\wedge\alpha_2'']) \\
  &  & +\lambda^2 (d\alpha_2'' + [\alpha_0\wedge\alpha_2'']),
\end{eqnarray*}
the coefficients of $\lm^{-1},\lm^0,\lm$ are respectively the projections of
$d\alpha + \alpha\wedge\alpha $ on $\g_{-1},\g_0,\g_1$ respectively
so they vanish and hence
$$
d\alpha_{\lambda} + \alpha_{\lambda}\wedge\alpha_{\lambda}=
d\beta_{\lambda^2} + \beta_{\lambda^2}\wedge\beta_{\lambda^2}
$$
where $\beta_\lm =\lm^{-1}\alpha_2' + \alpha_0 +\lm\alpha_2''$ is the
extended Maurer-Cartan form of $\beta=F^{-1}.dF$, the Maurer-Cartan
form of the lift $F\in G_I$  of $\rho_X\in S^{I}$. According to \cite{DPW},
we know that $\rho_X$ is harmonic \iif $d\beta_{\lambda} + \beta_{\lambda}\wedge
\beta_{\lambda}=0$, $\forall \lm\in \C^*$. This proves the second point
and  completes the proof.\hfill$\blacksquare$\\
\begin{rmk}\emph{
We have $\rho_-(x,y)=-\im(\overline x .y)=\rho_+(\overline x ,\overline y)$. Hence
$X\colon\Omega\to\h$ is  $\rho_-$-harmonic \iif $\overline X$ is
 $\rho_+ $-harmonic, and $X$ is $\omega_I^- $-isotropic \iif $\overline{X} $ is
$\omega_I^+$-isotropic. Besides if $U=(F,X)\colon\Omega\to
G_I\ltimes\h$ is a left LC$\omega_I$ then we have $F=L_a$ and
$aea^{-1}=\rho_X=\rho_+(q,q')$ with $dX=e^\omega(qdu +q'dv)$, $(q,q')$ orthonormal.
Thus $\rho_-(\overline q,\overline q')=aea^{-1}$ and hence
$\underline{U}=(R_{\overline a},\overline X)$ is a right LC$\omega_I$.}
\end{rmk}

\begin{rmk}\emph{
The restriction to $\im\h=\R^3$ of the left (or right) cross product gives us the usual
cross product in $\R^3$. Hence a surface in  $\im\h$ is left (resp.
right) $\rho$-harmonic \iif it is a constant mean curvature surface.\\
In the same way, it is easy to see that a surface in  $S^3$ is left (resp. right) $\rho$-harmonic \iif it is a constant mean curvature surface.}
\end{rmk}

\begin{rmk}\em
We can apply now the Dorfmeister-Pedit-Wu method (DPW) to obtain a Weierstrass representation
of $\rho$-harmonic surfaces (see \cite{DPW,HR1,HR3,ki1,ki2}). There are non-trivial technical difficulties in establishing DPW, such as proving loop group splittings (\cite{DPW,PS}).  
\end{rmk}

%%%%%%%%%%%%%%%%%%%%%%%%%%%%%%%%%%%%%%%%%%%%%%%%%%%%%%%%%%%%%%%%%%%%%%%%%%%%%
%                                                                           %
%       SECTION  : Second Elliptic Integrable System                        %                 %
%                                                                           %
%%%%%%%%%%%%%%%%%%%%%%%%%%%%%%%%%%%%%%%%%%%%%%%%%%%%%%%%%%%%%%%%%%%%%%%%%%%%%

\section{Second Elliptic Integrable Systems}\label{heart}
\subsection{4-symmetric spaces and twistor spaces}\label{twistor space}

\begin{defn}
Let $M$ be a Riemannian symmetric space. We will say that a Lie
group $G$  acts symmetrically on $M$ or that
$M$ is a $G$-symmetric space if $G$ acts
transitively and isometrically on $M$ and if there exists an involutive automorphism of
$G$, $\sigma$, such that $H$ the isotropy subgroup at a fixed point
$p_0\in M$, satisfies $(G^{\sigma})^0\subset H\subset
G^{\sigma}$. We will say also that $G/H$ is a symmetric realisation of
$M$.\\
We will say that a $G$-homogeneous space $N=G/G_0$ is a 4-symmetric bundle over the $G$-symmetric space $M$ if there exists an order four automorphism $\tau$ of $G$, such that $(G^{\tau})^0\subset G_0\subset G^{\tau}$, and
$(G,\tau)$ gives rise to the symmetric space $M$,  i.e. 
$\sigma=\tau^2$ and $G_0\subset H$.\\
A $G$-homogeneous space $N=G/G_0$ is a locally 4-symmetric space if there exists an order four 
automorphism of the Lie algebra 
$\g=\mathrm{Lie}\,G$, $\tau\colon\g\to\g$ such that $\g^{\tau}=\mathrm{Lie}\,G_0$. We will
say that $G/G_0$ is a locally 4-symmetric bundle over the $G$-symmetric
space $M$ if $\tau^2=\sigma$ (and $G_0\subset H$).
\end{defn}
Let us consider $M$ a $G$-symmetric space with
$\tau\colon\g\to\g$ an order four automorphism such that $\tau^2=\sigma$.
The automorphism
$\tau$ gives us an eigenspace decomposition of $\g^{\C}$:
$$\g^{\C}=\bigoplus_{k\in\mathbb{Z}_4}\gt_k$$
where $\gt_k$ is the $e^{ik\pi/2}$-eigenspace of $\tau$. We have
clearly $\gt_0=\g_0^{\C}$, $\overline{\gt_k}=\gt_{-k}$ and
$[\gt_k,\gt_l]\subset\gt_{k+l}$. We define $\g_2$, $\mathfrak{m}$ and $\underline{\g}_1$
 by
$$
\gt_2=\g_2^{\C}, \
\mathfrak{m}^{\C}=\gt_{-1}\oplus\gt_1 \,\text{ and }\,
\underline{\g}_1^{\C}=\bigoplus_{k\in\mathbb{Z}_4\smallsetminus\{0\}}\gt_k ,$$
it is possible because $ \overline{\gt_2}=\gt_2$ and $\overline{\gt_{-1}}
=\gt_1$. Let us set $\g_{-1}=\gt_{-1}$,\,$\g_1=\gt_1$ (i.e. we forget the " $\tl{}$ "), $\mathfrak{h}=
\g_0\oplus \g_2$ . Then
$$
\g=\mathfrak{h}\oplus\mathfrak{m}
$$
is the eigenspace decomposition of the involutive automorphism
$\sigma$, $\mathfrak{h}$ is the Lie algebra of $H$, the isotropy subgroup of $G$ at a reference point $p_0$, and $\mathfrak{m}$ is identified to  the tangent space $T_{p_0}M$. Besides we remark that $\tm\in\Sigma(\mk)$ (since $\tau_{|\mk^\C}=-i\mathrm{Id}_{\g_{-1}}\oplus i\mathrm{Id}_{\g_{1}}$)\footnote{We choose a metric in $\mk$ invariant by $\tm$ (and of course by $\Ad H$), see section~\ref{4-sym}}, which gives us the following theorem (proved in section~\ref{twistor}).
\begin{thm}\label{immersion}
Let us consider $M$ a Riemannian $G$-symmetric space and
$\tau\colon\g\to\g$ an order four automorphism such that $\tau^2=\sigma$. Let us make $G$ acting on $\Sigma(M)$: $g\cdot J=gJg^{-1}$. Let $J_0\in\Sigma(T_{p_0}M)$ be the 
complex structure corresponding\footnote{About the choice of $-\tm$ (instead of $\tm$) and its link to  the $(1,0)$-splitting, see theorem~\ref{fondamental} and remark~\ref{explanation}, for later explanation.} to $-\tm\in\Sigma(\mk)$, under the identification $T_{p_0}M=\mk$. Then the orbit of $J_0$ under the 
action of $G$ is an immersed submanifold of $\Sigma(M)$. Denoting by $G_0$ the stabilizer of $J_0$, then $\lie G_0=\g^\tau$ and thus $G/G_0$ is a locally 4-symmetric bundle over $M$, 
and the natural map
$$
 \begin{array}{crcl}
i \colon & G/G_0 & \longrightarrow &\Sigma(M)\\
 & g.G_0 & \longmapsto & gJ_0g^{-1}
  \end{array}
$$
is an injective immersion and a morphism of bundle. Moreover, if  the image of $G$ in $\Is(M)$ (the group of isometry of $M$) is closed, then $i$ is an embedding.
\end{thm}
%
%%%%%%%%%%%%%%%%%%%%%%%%%%%%%%%%%%%%%%%%%%%%%%%%%%%%%%%%%%%%%%%%%%%%%%%%%%%%%%%%%%%%%%%%%%%%%%%%%%%%%%%%%
%
\subsection{The second elliptic integrable system associated to a
4-symmetric space }\label{2.8}

We give ourself $M$ a Riemannian $G$-symmetric space with $\tau\colon\g\to\g$ an order four automorphism such that $\tau^2=\sigma$, and $N=G/G_0$ the associated locally 4-symmetric space 
 given by theorem~\ref{immersion}. We use the same notations as in  Section~\ref{twistor space}.
Then let us recall what is a second elliptic system according to C.L. Terng (see \cite{tern}).
\begin{defn}
The second $(\g,\tau)$-system is the equation for
$(u_0,u_1,u_2)\colon \C\to \oplus_{j=0}^2 \gt_{-j}$,
\begin{equation}\label{syst}
  \left\{
  \begin{array}{lr}
   \partial_{\bar{z}}u_2 + [\bar{u}_0,u_2]=0  & (a)\\
   \partial_{\bar{z}}u_1 + [\bar{u}_0,u_1] + [\bar{u}_1,u_2]=0  &  (b)\\
   -\partial_{\bar{z}}u_0 + \partial_{z}\bar{u}_0 + [u_0,\bar{u}_0] +
[u_1,\bar{u}_1] + [u_2,\bar{u}_2]=0.   & (c)
  \end{array}\right.
\end{equation}
It is equivalent to say that the 1-form
\begin{equation}\label{form}
\alpha_{\lm} =\sum_{i=0}^2 \lm^{-i} u_i dz + \lm^{i} \bar{u}_i
d\bar{z}=\lm^{-2}\alpha_2' + \lm^{-1}\alpha_1' + \alpha_0 +
\lm\alpha_1'' + \lm^2\alpha_2''
\end{equation}
satisfies the zero curvature equation:
\begin{equation}\label{courbnul}
d\alpha_{\lm} + \frac{1}{2}[\alpha_{\lm}\wedge\alpha_{\lm}]=0,
\end{equation}
for all $\lm\in \C^*$. We will speak about the $(G,\tau)$-system ($\tau$ is an
automorphism of $\lie G=\g$)  when we will look for solutions of the $(\g,\tau)$-system in $G$,
i.e. maps $U\colon\C\to G$ such that their Maurer-Cartan form is solution of the
$(\g,\tau)$-system, in other words when we integrate the zero curvature
equation (\ref{courbnul}) in $G$. We will call (geometric) solution of the second
elliptic integrable system associated to the locally 4-symmetric space $G/G_0$
a map $J\colon\C\to G/G_0$ which can be lifted into a solution $U\colon\C\to G$
of (\ref{syst}).
\end{defn}
\begin{rmk}\emph{
In (\ref{syst}), $\{\mrm{Im}((a)),(b),(c)\}$ is equivalent to $d\alpha +
\frac{1}{2}[\alpha\wedge\alpha]=0$. Hence the additional condition
added to the Maurer-Cartan equation by the zero curvature equation
(\ref{courbnul}) is  $\mrm{Re}\left(\partial_{\bar z}\alpha_2'( \dl{}{z}) 
+\left[\alpha_0''( \dl{}{\bar z}),\alpha_2'( \dl{}{z})\right]\right)=0$ or equivalently
$$
d(\star\alpha_2) + [\alpha_0\wedge(\star\alpha_2)]=0.
$$
}\end{rmk}
The first example of second elliptic system was given by F. H\'{e}lein
and P. Romon (see \cite{{HR1},{HR3}}): they showed that the
equations for Hamiltonian stationary Lagrangian surfaces in 4-dimension Hermitian
symmetric spaces are exactly the second elliptic system associated to
certain 4-symmetric spaces. Then in \cite{ki1}, we found another example in $\oct$:
the $\rho$-harmonic surfaces in $\oct$, which by restriction to $\h$ gave us the $\rho$-harmonic surfaces
in $\h$ (studied in section~\ref{r4}) which generalize the Hamiltonian stationary Lagrangian surfaces
in $\C^2$. 
\begin{defn}
Let $M$ be a Riemannian manifold and $\nabla$ its Levi-Civita
connection which induces a connection on $\mathrm{End}(TM)$. Let us
define for each $(p,J_p)\in\Sigma(M)$ the orthogonal projection
$$
\mathrm{pr}^{\bot}(p,J_p)\colon \mathrm{End}(T_pM)\to T_{J_p}(\Sigma(T_pM))
$$
($T_pM$ is an Euclidean vector space so $\Sigma(T_pM)$ is a
submanifold of the Euclidean space $\mathrm{End}(T_pM)$ and so
$T_{J_p}\Sigma(T_pM)$ is a vector subspace of $\mathrm{End}(T_pM)$
and we can consider the orthogonal projection on this subspace).
Given $L$ a Riemannian surface and $J\colon L\to\Sigma(M)$ we set
$$
\Delta J=\mathrm{pr}^{\bot}(J).\mrm{Tr}(\nabla^2 J)
$$
where $\mrm{Tr}$ is the trace with respect to the metric on $L$ (in fact, we take the vertical part of the rough Laplacian) .
We will say that $J$ is vertically harmonic if $\Delta J=0$. This notion depends only on the conformal structure on $L$. 
\end{defn}
\begin{defn}
Let $(L,j)$ be a Riemann surface, $M$ an oriented manifold and $X\colon
L\to M$ an immersion. Let $J\colon L\to X^*(\Sigma(M))$ be an
almost complex structure on the vector bundle $X^*(TM)$. Then we
will say that  $J$ is an admissible twistor lift of $X$ if one of the
following equivalent statements holds:
\begin{description}
  \item[(i)] $X$ is $J$-holomorphic: $\star dX:=dX\circ j=J.dX$
  \item[(ii)] $J$ is an extension of the complex structure on the
  oriented tangent plane $P=X_{*}(TL)$ induced by $j$, the complex
  structure of $L$, or equivalently $J$ induces the complex
  structure $j$ in $L$.
  \item[(iii)] $X$ is a conformal immersion and $J$ stabilizes the
  tangent plane $X_{*}(TL)$, i.e. for all $z\in L$, $J_z$
  stabilizes $X_{*}(T_zL)$ and induces on it the same
  orientation, which we will denote by $J\circlearrowleft X_{*}(TL)$
  \item[(iv)] $X$ is a conformal immersion and $J$ is an extension
  of the unique positive complex structure $I_P$ of the tangent plan
  $P=X_{*}(TL)$.
\end{description}
Finally, we will say that a map $J\colon L\to\Sigma(M)$ is an
admissible twistor lift if its projection $X=\mrm{pr}_M\circ J\colon
L\to M$ is an immersion and $J$ is an admissible twistor lift of it.
\end{defn}

\begin{thm}\label{fondamental}
Let $L$ be a simply connected Riemann surface and $(G,\tau)$
a locally 4-symmetric bundle over a symmetric space $M=G/H$. Let
$J_0\in\Sigma(T_{p_0}M)$ be the complex structure corresponding to
$-\tau_{|\mk}$ (see Section \ref{twistor space}). Let be $J_X\colon L\to i(G/G_0)\subset
\Sigma(G/H)$.
 Then the two following statements are equivalent:
\begin{description}
  \item[$\bullet$] $J_X$ is  an admissible twistor lift.
  \item[$\bullet$] Any lift $F\colon L\to G$ of $J_X$
  ($FJ_0 F^{-1}=J_X$) gives rise to a
  Maurer-Cartan form $\alpha=F^{-1}.dF$ which satisfies:
  $\alpha_{-1}''=\alpha_1'=0$ and $\alpha_{-1}'$ does not vanish.
  \end{description}
Furthermore, under these statements, $J_X\colon L\to\Sigma(G/H)$ is vertically harmonic
\iif $J_X\colon L\to G/G_0$ is solution of the second elliptic
integrable system associated to the locally 4-symmetric space
$(G,\tau)$, i.e.
$$
d\alpha_{\lm} + \frac{1}{2}[\alpha_{\lm}\wedge\alpha_{\lm}]=0,
\quad\forall\lm\in \C^*,
$$
where $\alpha_{\lm}=\lm^{-2}\alpha_2' + \lm^{-1}\alpha_{-1}' + \alpha_0 +
\lm\alpha_1'' + \lm^2\alpha_2''$
is the extended Maurer-Cartan form of $\alpha$.
\end{thm}
\textbf{Proof.} For the first point, let us  make $F^{-1}$ acting on
the equation $dX\circ j=J_X.dX$, we obtain
$\alpha_{\mk}\circ j=-\tm(\alpha_{\mk})$ which is equivalent to
$\alpha_{-1}''=\alpha_1'=0$. Thus
$\alpha_{-1}(\dl{}{z})=\alpha_{\mk}(\dl{}{z})=F^{-1}.\dl{X}{z}$, and $X$ is an
immersion \iif $\alpha_{-1}'$ does not vanish.\\
For the second point, let us recall that
$\mrm{End}(T_{p}M)=\mrm{sym}(T_pM)\overset{\bot}{\oplus}\so(T_pM)$
and given $J\in\Sigma(T_pM)$, we have
$T_J{\Sigma(T_pM)}=\mrm{Ant}(J)=\{A\in\so(T_pM)| AJ+JA=0\}$ and $(T_J{\Sigma(T_pM)}
)^{\bot}\cap\so(T_pM)=\mrm{Com}(J)=\{A\in\so(T_pM)|[A,J]=0\}$.\\
Now, let us compute the connection $X^*\nabla$ on
$X^*(\mrm{End}(TM))$, in terms of the Lie algebra setting. Let $A$ be a
section of $X^*(\mrm{End}(TM))$ and $Y$ a section of $X^*(TM)$. Let
$A_0\in C^\infty(L,\mrm{End}(T_{p_0}M))$ be  defined by
$A_{F.p_0}=FA_0F^{-1}$ and $A_{\mk}\in C^\infty(L,\mrm{End}(\mk))$ its
image under the identification $T_{p_0}M=\mk$. Then $A_{F.p_0}$
corresponds to $\Ad F\circ A_{\mk}\circ\Ad F^{-1}$ (under the
identification $TM=[\mk]:=\{(g.p_0,\Ad g(\xi)),\, \xi\in\mk ,g\in
G\}$, see section~\ref{4-sym}). In particular $(J_X)_{\mk}=-\tm$ (we
suppose $F(p_0)=1$). We set also $Y=\Ad F(\xi).p_0$, $\xi\in
C^{\infty}(L,\mk)$. From now, we do the identification
$TM=[\mk]$ without precising it. Then, denoting  by $[\ ,\ ]_\mk$ the $\mk$-component of the Lie bracket, we have
\begin{eqnarray*}
(\nabla A)(Y) & = & \nabla(AY) -A(\nabla Y) \\
              & = & \Ad F\left([d(A_{\mk}\xi) + [\alpha,A_{\mk}.\xi]]_{\mk}
 - A_{\mk}(d\xi + [\alpha,\xi]_{\mk})\right)\\
              & = & \Ad F \left((dA_{\mk})\xi +
 (\ad\alpha_{\hk}\circ A_{\mk}- A_{\mk}\circ\ad\alpha_{\hk})\xi\right).
\end{eqnarray*}
Hence
$$
\nabla A= \Ad F(dA_{\mk} + [\ad_{\mk}\alpha_{\hk},A_{\mk}]).
$$
In particular,\footnote{In all the proof, we will merge $\alpha_k'$ (resp. $\alpha_k''$) with $\alpha_k'(\dl{}{z})$ (resp. 
$\alpha_k''(\dl{}{\bar z})$), and in particular write `$[\alpha_k'',\alpha_l']$' instead of `$[\alpha_k''(\dl{}{\bar z}),
\alpha_l'(\dl{}{z})]$'. $z$ is a local holomorphic coordinate in $L$.}
$$\nabla_{\dl{}{z}}J_X=-2\Ad F(\ad_{\mk}\alpha_2'\circ\tm)
$$
(because $\ad_{\mk}\g_0$ commutes with
$\tm$ whereas $\adm\g_2$ anticommutes with it) and thus
\begin{eqnarray*}
\nabla_{\dl{}{\bar z}}(\nabla_{\dl{}{z}}J_X) & = & -2\Ad F\left(\ad_{\mk}(\partial_{\bar z}
\alpha_2')\circ\tm  + [\ad_{\mk}(\alpha_{\hk}''),\ad_{\mk}(\alpha_2')\circ\tm]\right)\\
     & = &     -2\Ad F\left(\ad_{\mk}(\partial_{\bar z}\alpha_2')\circ\tm
+ \ad_{\mk}([\alpha_0'',\alpha_2'])\circ\tm \right.\\
 &  &  \qquad\qquad\qquad\qquad\qquad +\left. [\ad_{\mk}\alpha_2'',\ad_{\mk}(\alpha_2')
 \circ\tm ]\right )\\
  & = &  -2\Ad F\left(\ad_{\mk}(\partial_{\bar z}\alpha_2'
+ [\alpha_0'',\alpha_2'])\circ\tm + [\ad_{\mk}\alpha_2'',\ad_{\mk}(\alpha_2')
\circ\tm ]\right )
\end{eqnarray*}
but $-\Ad F\left([\ad_{\mk}\alpha_2'',\ad_{\mk}(\alpha_2')\circ\tm
]\right)$ commutes with $-\Ad F(\tm)=J_X$ so it is orthogonal to
$T_J\Sigma(T_pM)$ thus
$$
\mrm{pr}^{\bot}(J_X).\nabla_{\dl{}{\bar z}}(\nabla_{\dl{}{z}}J_X) =-2\Ad F\left(\ad_{
\mk}(\partial_{\bar z}\alpha_2' + [\alpha_0'',\alpha_2'])\circ\tm\right ).
$$
Hence, since $\ad_{\mk}$ is injective\footnote{We can do this hypothesis without loss of generality, see section~\ref{4-sym}.}
\begin{equation}\label{equiv}
\Delta J_X=0\Longleftrightarrow \mrm{Re}\left(\partial_{\bar z}\alpha_2'
+ [\alpha_0'',\alpha_2']\right)=0.
\end{equation}
This completes the proof.\hfill$\blacksquare$

\begin{rmk}\emph{
The equivalence (\ref{equiv}) holds for any map $J_X\colon L\to
i(G/G_0)$. Indeed, we have not used the fact that $J_X$ is an
admissible twistor lift to prove this equivalence.}
\end{rmk}

\begin{thm} Let $J_X\colon L\to G/G_0\hookrightarrow\Sigma(G/H)$ be an admissible twistor
lift. Then $J_X\colon L\to G/G_0$ is harmonic\footnote{with respect to any metric induced by an $\Ad G$-invariant metric in $\g$.} \iif $X\colon L\to G/H$ is harmonic and
$J_X$ is vertically harmonic.
\end{thm}
\textbf{Proof.} $J_X\colon L\to G/G_0$ is harmonic \iif the Maurer-Cartan form $\alpha=F^{-1}.dF$ of the lift $F\colon L\to G$ of $J_X$
($FJ_0 F^{-1}=J_X$) satisfies (see \cite{12})
$$\displaystyle
{\partial_{\bar z}\underline{\alpha}_1' +
[\alpha_0'',\underline{\alpha}_1'] +
\frac{1}{2}[\underline{\alpha}_1'',\underline{\alpha}_1']_{\underline{\g}_1}=0}
$$
(where $\g=\g_0\oplus\underline{\g}_1$ is the reductive decomposition corresponding
to the homogeneous space $G/G_0$, see Section~\ref{twistor space}) which
splits into
\begin{equation}
\displaystyle
\left\{
  \begin{array}{l}
\displaystyle   \partial_{\bar{z}}\alpha_2' + [\alpha_0'',\alpha_2']+ \frac{1}{2}[\alpha_{1}'',
     \alpha_{1}']+
  \frac{1}{2} [\alpha_{-1}'',\alpha_{-1}'] =0 \\
\displaystyle   \partial_{\bar{z}}\alpha_{-1}' + [\alpha_0'',\alpha_{-1}'] +
\frac{1}{2}[\alpha_2'',\alpha_{1}']   + \frac{1}{2}[\alpha_{1}'',\alpha_2']=0  \\
\displaystyle   \partial_{\bar{z}}\alpha_{1}' + [\alpha_0'',\alpha_{1}'] +
  \frac{1}{2}[\alpha_2'',\alpha_{-1}']+ \frac{1}{2}[\alpha_{-1}'',\alpha_2']=0.
  \end{array}\right.
\end{equation}
then, using $\alpha_{-1}''=\alpha_1'=0$, we obtain
$$
\displaystyle
\left\{
  \begin{array}{l}
\partial_{\bar z}\alpha_2'+ [\alpha_0'',\alpha_2']  =  0 \\
 \partial_{\bar{z}}\alpha_{-1}' + [\alpha_0'',\alpha_{-1}']  =  0 \\
     {[\alpha_{2}'',\alpha_{-1}']}  =  0
\end{array}\right.
$$
(in the second equation, we have used $[\alpha_{1}'',\alpha_2']=-\overline{[
\alpha_{2}'',\alpha_{-1}']}=0$).\\
Besides $X\colon L\to G/H$ is harmonic \iif we have
$$
\partial_{\bar z}\alpha_{\mk}' +[\alpha_{\hk}'',\alpha_{\mk}']=0
$$
which splits into
\begin{equation}
\displaystyle
\left\{
  \begin{array}{l}
\displaystyle   \partial_{\bar{z}}\alpha_{-1}' + [\alpha_0'',\alpha_{-1}'] +
[\alpha_2'',\alpha_{1}']  =0  \\
\displaystyle   \partial_{\bar{z}}\alpha_{1}' + [\alpha_0'',\alpha_{1}'] +
   [\alpha_2'',\alpha_{-1}']=0.
  \end{array}\right.
\end{equation}
and using $\alpha_{-1}''=\alpha_1'=0$, we obtain
$$
\left\{
  \begin{array}{l}
\displaystyle   \partial_{\bar{z}}\alpha_{-1}' + [\alpha_0'',\alpha_{-1}']=0  \\
\displaystyle  [\alpha_2'',\alpha_{-1}']=0.
  \end{array}\right.
  $$
 This completes the proof.\hfill$\blacksquare$

%%%%%%%%%%%%%%%%%%%%%%%%%%%%%%%%%%%%%%%%%%%%%%%%%%%%%%%%%%%%%%%%%%%%%%%%%%%%%%%%%%%%%%%%%%%%
%                                                                                          %
%          Structure of 4-symmetric bundles over symmetric spaces                          %
%                                                                                          %
%                                                                                          %
%%%%%%%%%%%%%%%%%%%%%%%%%%%%%%%%%%%%%%%%%%%%%%%%%%%%%%%%%%%%%%%%%%%%%%%%%%%%%%%%%%%%%%%%%%%%
%
%
%
\section{Structure of 4-symmetric bundles over symmetric spaces}\label{structure}
\subsection{4-symmetric spaces}\label{4-sym}
Let $G$ be a Lie group with Lie algebra $\g$,
$\tau\colon G\to G$ an order four automorphism with the fixed point
subgroup $G^{\tau}$, and the corresponding Lie algebra
$\g_0=\g^{\tau}$. Let  $G_0$ be a subgroup of $G$ such that $(G^{\tau})^0\subset
G_0\subset G^{\tau}$, then $\mathrm{Lie}\,G_0=\g_0$ and $G/G_0$ is a 4-symmetric space. The automorphism
$\tau$ gives us an eigenspace decomposition of $\g^{\C}$ for which we use the notation of section~\ref{twistor space}.
Then $\g=\mathfrak{h}\oplus\mathfrak{m}$ is the eigenspace decomposition of the involutive automorphism
$\sigma=\tau^2$. Let
 $H$ be a subgroup of $G$ such that $(G^{\sigma})^0\subset H\subset G^{\sigma}$
then $\mathrm{Lie}\,H=\mathfrak{h}$ and $G/H$ is a symmetric space.
We will often suppose that $G_0$ and $H$ are chosen such that
$G_0=G^{\tau}\cap H$. With this condition, $G_0\subset H$ so that
$G/G_0$ is a bundle over $G/H$.
Recall that the tangent bundle  $TM$ is canonically isomorphic to
the subbundle $[\mk]$ of the trivial bundle $M\times\g$, with fibre
$\text{Ad}g(\mk)$ over the point $x=g.H\in M$. Under this
identification the canonical $G$-invariant connection of $M$ is just the flat
differentiation in $M\times\g$ followed by the projection on
$[\mk]$ along $[\hk]$ (which is defined  in the same way as
$\mk$) (see  \cite{8}).
For the homogeneous space $N=G/G_0$ we have the following reductive decomposition
\begin{equation}\label{reductive}
\g=\g_0\oplus\underline{\g}_1
\end{equation}
($\underline{\g}_1$ can be written $\underline{\g}_1=\mk\oplus\g_2$) with
$[\g_0,\underline{\g}_1]\subset\underline{\g}_1$. As for the symmetric space $G/H$,
we can identify the tangent bundle $TN$ with the subbundle $[\underline{\g}_1]$ of the
trivial bundle $N\times\g$, with fibre $\text{Ad}g(\underline{\g}_1)$ over the
point $y=g.G_0\in N$.\\
The symmetric space $M=G/H$ is Riemannian if it admits a $G$-invariant metric, which is
equivalent to say that $\mk$  admits an $\mathrm{Ad}(H)$-invariant 
inner product or equivalently, that $\mathrm{Ad}_{\mk}(H)$  be relatively compact\footnote{In the literature,
it is often supposed that $\mathrm{Ad}_{\mk}(H)$ is compact. We will see that these two hypothesis are in fact
equivalent.}. We remark that the Levi-Civita connection coincides with the previous canonical $G$-invariant connection and in particular is independent of the $G$-invariant
metric chosen. We will always suppose that  the symmetric  spaces $M$ which we
consider are Riemannian. We will in addition to that suppose that the $\Ad (H)$-invariant inner product in $\mk$ is also invariant by $\tm$ (such an inner product always exists when $\mathrm{Ad}_{\mk}(H)$ is relatively compact, see the appendix).
We will also suppose that $M$ is connected, then $G^0$ acts transitively on $M$ and so we can
suppose that $G$ is connected.\\
We want to study the Riemannian symmetric spaces $M$ such that there
exists a 4-symmetric space $(G,\tau)$ which gives rise to $M$ in the
same way as above. For that, let us recall the following theorem:
\begin{thm}\cite{besse,Hel}
Let $M$ be a Riemannian manifold.
\begin{description}
  \item[(a)] The group $\Is(M)$ of all the isometries of $M$ is a
  Lie group and acts differentiably on $M$.
  \item[(b)] Let $p_0\in M$, then an isometry $f$ of $M$ is
  determined by the image $f(p_0)$ of the point $p_0$ and the
  corresponding tangent map $T_{p_0}f$ (i.e. if $f(p_0)=g(p_0)$ and
  $T_{p_0}f=T_{p_0}g$ then $f=g$).
  \item[(c)] The isotropy subgroup $\Is_{p_0}(M)=\{f\in\Is(M);f(p_0)=p_0\}$
is a closed subgroup of $\Is(M)$ and the linear isotropy
representation $\rho_{p_0}\colon f\in \Is_{p_0}(M)\mapsto T_{p_0}f\in
O(T_{p_0}M)$ is an isomorphism from $\Is_{p_0}(M)$ onto a closed
subgroup of $O(T_{p_0}M)$. Hence $\Is_{p_0}(M)$ is a compact
subgroup of $\Is(M)$.
  \item[(d)]If $M$ is a Riemannian homogeneous space, $M=G/H$ with
  $G=\Is(M)$, $H=\Is_{p_0}(M)$ and $\mk$ an $\Ad H$-invariant space
  such that $\g=\hk\oplus\mk$, then the previous closed subgroup,
  image of $H$ by the preceding isomorphism $\rho_{p_0}$, i.e. the
  linear isotropy subgroup $H^*$ can be identified to $\Ad_{\mk}H$.
  More precisely the linear isometry $\xi\in \mk \mapsto\xi.p_0\in
  T_{p_0}M$ gives rise to an isomorphism from $O(\mk)$ onto
  $O(T_{p_0}M)$ which sends $\Ad_{\mk}H$ onto $H^*$. Hence the
  linear adjoint representation of $H$ on $\mk$: $g\in H\mapsto
  \Ad_{\mk}g\in\Ad_{\mk}H$ is an isomorphism (of Lie groups).
  $H\cong H^*\cong \Ad_{\mk}H$.
\end{description}
\end{thm}

\subsubsection{First convenient hypothesis.}\label{first}
There may be more than one Lie group $G$ acting symmetrically on a Riemannian
symmetric space $M$.
Besides, we have a convenient way to work on Riemannian symmetric
spaces: it is to consider  that $G$ is a
subgroup of the group of isometries of $M$, $\Is(M)$, which is
equivalent to suppose that $G$ acts effectively on $M$, i.e. $H$, the
isotropy subgroup at a fixed point $p_0$ does
not  contain non-trivial normal subgroup of $G$ (see \cite{besse}).
It is always possible because the kernel $K$ of the natural morphism $\phi_{H}\colon G\to
\Is(M)$ is the maximal normal subgroup of $G$ contained in $H$ \footnote{$K=\ker\phi_{H}
=\ker\rho_{p_0}=\ker\Ad_{\mk}$}, and $G'=G/K$ acts
transitively and effectively on $M=G/H$ with isotropy subgroup
$H'=H/K$. Thus $M=G'/H'$ and since $K\subset H\subset G^\sigma$, then $\sigma$ gives rise to an involutive
morphism $\sigma'\colon G'\to G'$ such that $({G'}^{\sigma'})^0\subset
H'\subset {G'}^{\sigma'}$. Now, let us suppose that there exists an
 order four automorphism $\tau\colon G\to G$ such that $\sigma=\tau^2$. Then it gives rise to
an isomorphism $\tau'\colon G/K \to G/\tau(K)$. We would like that
$\tau(K)=K$. It is the case if $\tau(H)=H$: $K$ and
$\tau(K)$ are respectively the maximal normal subgroups of $G$ contained in $H$ and
$\tau(H)$ respectively, and so if $\tau(H)=H$ then $K=\tau(K)$.\\[1.5mm]
Let us suppose that $\tau(K)=K$, then $\tau$ gives rise to
an order four automorphism $\tau'\colon G/K\to G/K$ such that
$\sigma'={\tau'}^2$. 
With our convention we have $G_0'={G'}^{\tau'}\cap  H'$,
then we obtain a 4-symmetric bundle $N_{min}'=G'/G_0'$ over $M$.
Hence, when $G_0'$ describes all the possible choices:
$({G'}^{\tau'})^0\subset G_0'\subset {G'}^{\tau'}\cap H'$, we obtain
a family of 4-symmetric bundles $N'=G'/G_0'$ over $M$ which are discrete coverings
of $N_{min}'=G'/({G'}^{\tau'}\cap H')$ and of which $N_{max}'=G'/({G'}^{\tau'})^0$ is a
discrete covering. For example, if we choose $G_0'=\pi_{K}(G_0K)$, we
obtain the 4-symmetric bundle over $M$, $N'=(G/K)/\pi_{K}(G_0K)=G/G_0K=N/K$ \footnote{In the writing $N'=N/K$, $K$ does not act freely on $N$ in general: it is 
$K'=K/K\cap G_0$ which acts freely on $N$ and we have $N'=N/K=N/K'$. In particular it is possible that $N/K=N$ for a non-trivial $K$ (see section~\ref{cgrass}).}.\\[1.5mm]
\indent Let us come back to the general case (i.e. we do not suppose that
$\tau(K)=K$).\\[0.8mm]
Since $\tau(\hk)=\hk$, we have $\tau(H^0)=H^0$ and thus denoting by $K_0$ the
maximal normal subgroup of $G$ contained in $H^0$ (we have $K^0\subset K_0\subset
K\cap H^0$), then $\tau(K_0)=K_0$ for the same reason as above (in
particular, if $K_0=K$ i.e. $K\subset H^0$, then we are in the
previous case: $\tau(K)=K$).
Hence $\tau$ gives rise to an order four automorphism
$\tilde{\tau}\colon G/K_0\to G/K_0$ and we are in the case considered above
if we consider the symmetric space $\tilde{M}=G/H^0$ (instead of $M$). Let us precise
this point. Indeed $\tilde{M}$ is  a $(G/K_0)$-symmetric space and $\tilde{G}=G/K_0$ acts
effectively on it (the isotropy subgroup $\tilde{H}=H^0/K_0$ does not contain non-trivial normal
subgroup of $G/K_0$): as above $\sigma$ gives rise to an involutive automorphism
$\tilde{\sigma}$  of $\tilde{G}=G/K_0$ such that $\tilde{H}=(\tilde{G}^{\tilde{
\sigma}})^0$  and $\tilde{\tau}$ is
an order four automorphism of $G/K_0$ such that $\tilde{\tau}^2=\tilde{\sigma}$.
Finally, as above we obtain a family of 4-symmetric bundles
$\tilde{N}=\tilde{G}/\tilde{G}_0$ over $\tilde{M}$ when $\tilde{G}_0$
describes the set of all possible choices:
$({\tilde{G}}^{\tilde{\tau}})^0\subset \tilde{G}_0\subset
{\tilde{G}}^{\tilde{\tau}}\cap \tilde{H}$.\\
Moreover, 
the involution  $\tilde{\sigma}$ of $G/K_0$  gives
rise also to the $G/K_0$-symmetric space $M$ (i.e. $({\tilde{G}}^{\tilde{\sigma}})^0\subset
H/K_0\subset {\tilde{G}}^{\tilde{\sigma}}$  or equivalently $M$ belongs to the
family of  $G/K_0$-symmetric spaces defined by $\tilde{\sigma}$ (of
which $\tilde{M}$ is a discrete covering)).\\[1mm]
In the same way, we have $\tau(G^{\sigma})=G^{\sigma}$ and thus we can
do the same as above for the symmetric space $M_{min}=G/G^{\sigma}$.\\[1.5mm]
\indent Nevertheless, in general, it is possible that $\tau(K)\neq K$ and then
$\tau$ does not give rise to an order four automorphism of $G'=G/K$ but only to the
isomorphism $\tau'\colon G/K \to G/\tau(K)$. However, the tangent map $T_e\tau'=
T_e\tilde{\tau}$ is an order four automorphism of the Lie algebra
$\mathrm{Lie} (G/K) = \mathrm{Lie} (G/\tau(K))=\mathrm{Lie}(G/K_0)=\g/\mathfrak{k}$,
and we have $(T_e\tau')^2=T_e\sigma'$, thus $N/K=(G/K)/\pi_K(G_0K)$
is a locally 4-symmetric bundle over $M$ ($\lie \pi_K(G_0K)=\g^{T_e\tau'}$).

Hence we have two good settings to study the Riemannian
symmetric spaces $M$ over which a 4-symmetric bundle can be defined, if we
want to work only with subgroups of $\Is(M)$.\\[1.5mm]
\textbf{The first possibility} is to consider that we begin by giving
ourself an order four automorphism $\tau\colon G\to G$ and that we
always choose the Riemannian symmetric space $\tilde{M}=G/H$ with
$H=(G^{\tau^2})^0$ (respectively $M_{min}=G/H$ with $H=G^{\tau^2}$). In other words, in the family of $G$-symmetric space
corresponding to $\sigma=\tau^2$ (i.e. $(G^{\sigma})^0\subset H\subset
G^{\sigma}$), we choose the "maximal" one $\tilde{M}=G/(G^{\sigma})^0$, which is a
discrete covering of all the others (respectively the "minimal" one $M_{min}=G/G^{\sigma}$,
of which all the others are discrete coverings). Then according to what precedes,
we can always suppose that $G$ is a subgroup of $\Is(\tilde{M})$ (respectively of $\Is(M_{min})$).\\[1.5mm]
\textbf{The second possibility} is to work with locally 4-symmetric spaces. In other words
we begin by a Riemannian symmetric space over which there exists a  locally 4-symmetric
bundle. It means that we work with the following setting: a Riemannian
symmetric spaces $M$ with $G$ a subgroup of $\Is(M)$ acting
symmetrically on $M$ and an order four automorphism $\tau\colon\g\to\g$,
such that $\tau^2=\sigma$. To define the locally 4-symmetric
space $N$ in this setting, we must tell how we define $G_0$. We will
set
\begin{equation}\label{g-0}
 G_0=\{g\in H |\,\Ad_{\mk} g\circ\tau_{|\mk}\circ\Ad_{\mk} g^{-1}=\tau_{|\mk}\}.
 \end{equation}
First, we have to verify that if $\tau$ can be integrated by an
automorphism of $G$, also denoted by $\tau$, then we have
$G_0=G^{\tau}\cap H$. Indeed, if $g\in G^{\tau}\cap H$, then $\Ad g\circ\tau\circ
\Ad g^{-1}=\Ad (g.\tau(g)^{-1})\circ\tau=\tau$ and since $\Ad H$ stabilizes
$\mk$, we have $\Ad_{\mk} g\circ\tau_{|\mk}\circ\Ad_{\mk} g^{-1}=\tau_{|\mk}$ by
taking the restriction to $\mk$ of the preceding equation.
Conversely, suppose that $g\in H$ and $\Ad_{\mk} g\circ\tau_{|\mk}\circ\Ad_{\mk} g^{-1}=
\tau_{|\mk}$, then $\Ad (g.\tau(g)^{-1})\circ\tau_{|\mk}=\tau_{|\mk}$ so
since $\tau_{|\mk}$ is surjective, $\Ad(g.\tau(g)^{-1})_{|\mk}=\mathrm{Id}_{\mk}$
and since the adjoint representation of $H$ on $\mk$ is injective
(because we suppose that $G$ is a subgroup of $\Is(M)$, and thus $H$ is a subgroup of
$\Is_{p_0}(M)$) it follows that $g.\tau(g)^{-1}=1$. Finally, $g\in
G^{\tau}\cap H$. Thus our definition (\ref{g-0}) is coherent with our convention
which holds when $\tau$ can be integrated by an automorphism of $G$.\\
Besides, it is easy to see that $\lie
G_0=\{a\in\hk|\,\adm a\circ\tau_{|\mk}=\tau_{|\mk}\circ\adm a\}=\g_0$.
(Indeed, $\forall a\in \g_0$, $\ad a\circ\tau=\tau\circ\ad a$,
and $\forall a\in \g_2$, $\ad a\circ\tau=-\tau\circ\ad a$,
moreover $\tau_{|\mk}\circ {\adm a}=0\Rightarrow {\adm a}=0\Rightarrow
a=0$ because $a\in\hk\mapsto {\adm a}$ is the tangent map of
$h\in H\mapsto \Adm h$ which is an injective morphism). Hence
$N=G/G_0$ is a locally 4-symmetric bundle over $M$.\\
Further, let $\pi\colon\tilde{G}\to G$ be the universal covering of $G$, and
$D=\ker\pi$. Then $\tau$ can be integrated by $\tilde{\tau}\colon\tilde{G}\to G$.
Set $\tilde{\sigma}=\tilde{\tau}^2$, then $\sigma\circ\pi=\pi\circ\tilde{\sigma}$
and $T_1\sigma=T_1\tilde{\sigma}=(T_1\tilde{\tau})^2$. $\tilde{G}$
acts almost effectively on $M$ with isotropy subgroup
$\tilde{H}=\pi^{-1}(H)$ and almost effectively on $\tilde{M}=\tilde{G}/\tilde{H}^0$
which is the universal covering of $M$ (see \cite{Hel}). Besides, if $\tilde{G}$ does
not act effectively on $\tilde{M}$, then we take $D_0$ the maximal
normal subgroup of $\tilde{G}$ included in $\tilde{H}^0$, and then we
quotient by it, so that we obtain an effective action of $\tilde{G}/D_0$ on $\tilde{M}$
and $\tilde{\tau}$ gives rise to an automorphism of $\tilde{G}/D_0$,
according to above.
Thus we are in the first possibility. Besides
it is easy to see that $\forall g\in\tilde{G}$, $\Ad g =\Ad\pi(g)$
(more precisely $T_1\pi\circ\Ad g=\Ad\pi(g)\circ T_1\pi$ and we
identify $\gt$ and $\g$ so that $T_1\pi=\mathrm{Id}$). Thus
$\tilde{G}_0=\tilde{G}^{\tilde{\tau}}\cap
\tilde{H}^0=\{g\in\tilde{H}^0|\,\Ad g\circ\tau_{|\mk}\circ\Ad
g^{-1}=\tau_{|\mk}\}\subset \pi^{-1}(G_0)$ . Hence the 4-symmetric
space $\tilde{G}/\tilde{G}_0$ is a discrete covering of the locally
4-symmetric space $G/G_0$ and we have the following commutative
diagram:
\begin{equation}\label{diag}
\begin{CD}
\tilde{G}/\tilde{G}_0 @>>> G/G_0 \\ @VVV
@VVV \\
\tilde{M} @>>> M\ .
\end{CD}
\end{equation}
\textbf{In conclusion}, the two possibilities are equivalent, but we will use the
second one because it works with any symmetric space $M$, whereas
the first one needs that we choose a certain covering of $M$ (for example
its universal covering).

\begin{rmk}\label{univcov}\emph{
We see that in the preceding reasoning (this using the universal covering
$\tilde{G}$) we need only the automorphism of Lie algebra $\tau$
(and not the symmetric space $M$). Hence, we can consider that we
work in the Lie algebra setting and give ourself an order four
automorphism $\tau$ of $\g$. Then we consider the family of
associated pairs $(G,H)$ where $G$ is a connected Lie group with Lie
algebra $\g$ and $H$ a closed Lie subgroup with Lie algebra
$\hk=\g^\sigma$. To each such pair
corresponds the locally symmetric space $M=G/H$ and defining $G_0$
by (\ref{g-0}), the locally 4-symmetric bundle $N=G/G_0$ over $M$.
Let $\tilde{G}$ be a simply connected Lie group with Lie algebra
$\g$, then $\tau$ and $\sigma$ integrate in $\tilde{G}$ and thus
for  $\tilde{H}$ the closed subgroup with Lie algebra $\hk$, we
can take any subgroup such that
$(\tilde{G}^{\tilde{\sigma}})^0\subset\tilde{H}\subset\tilde{G}^{\tilde{\sigma}}$
(which implies that $\tilde{H}$ is closed). If we suppose $\tilde{H}$
connected, i.e. $\tilde{H}=(\tilde{G}^{\tilde{\sigma}})^0$, then
$\tilde{M}=\tilde{G}/\tilde{H}$ is a symmetric space and is also the
universal covering of all the locally symmetric spaces $M=G/H$
when $(G,H)$ describes all the associated pairs (see \cite{Hel}),
and we have the above commutative diagram between the 4-symmetric bundle $\tilde{N}=\tilde{G}/
\tilde{G}_0$ over $\tilde{M}$ and the locally 4-symmetric bundle $N=G/G_0$ over
$M$. Moreover if $\tl M$ is Riemannian then all the symmetric spaces $M=G/H$ when $(G,H)$ describes all
the symmetric associated pairs are Riemannian (see appendix, corollary~\ref{cpair}).
}\end{rmk}

\begin{rmk}\label{tau_h}\emph{
Let us consider $M$ a $G$-symmetric space, $G\subset\Is(M)$, and $\tau\colon\g
\to\g$ an order four automorphism such that $\tau^2=\sigma$. Then
we have $\tau_{|\mk}\in\Sigma(\mk)$ ($\tau_{|\mk^\C}=-i\mathrm{Id}_{\g_{-1}}\oplus
i\mathrm{Id}_{\g_{1}}$) and it is easy to see that
$$
\forall a\in\hk,\
\tau_{|\hk}(a)=\ad_{\mk}^{-1}(\tau_{|\mk}\circ\adm a\circ\tau_{|\mk}^{-1}).
$$
In other words, under the identification  $\hk\simeq\ad_{\mk}\hk\subset \mathfrak{so}(\mk)$, $\tau_{|\hk}$ is
the restriction to $\hk$ of $\Ad (\tau_{|\mk})\colon \so(\mk)\to \so(\mk)$.
Hence $\tau$ is determined by $\tau_{|\mk}$. Besides $\tau_{|\hk}$ is the tangent
map  of the isomorphism $\tau_H$:
$$\tau_H(g)=\Ad_{\mk}^{-1}(\tau_{|\mk}\circ\Ad_{\mk} g\circ\tau_{|\mk}^{-1}),
$$
for $g\in H^0$ (and more generally for $g\in\Ad_{\mk}^{-1}\circ(\Int\tm)^{-1}\circ\Ad_{\mk}(H)$). Under the identification $H\simeq\Ad_{\mk}H\subset
O(\mk)$ it is the restriction to $H^0$ of the involution $\Int\tm\colon
O(\mk)\to O(\mk)$. According to the definition (\ref{g-0}) of $G_0$, we
have $G_0=H^{\tau_H}$. Besides $\tau_H(H^0)=H^0$, thus $H^0/G_0^0$ is an $H^0$-symmetric space. If $\Int\tm(\Ad_{\mk} H)=(\Ad_{\mk} H)$, then $
\tau_H$ is defined in $H$ and $\tau_H(H)=H$, then $H/G_0$ is an $H$-symmetric
space (if $\tau_H(H)\neq H$ it is only a locally symmetric space). Obviously, if $\tau$
can be integrated in $G$ then $\tau_H=\tau_{|H}$.}\end{rmk}
\begin{defn}\label{autm}
Let $M$ be a $G$-symmetric space.
Let  $\mrm{Aut}(\mk)$ be the subgroup of $O(\mk)$ defined by:
$$
\mrm{Aut}(\mk)=\{F\in O(\mk)\mid F(\ad_{\mk}[v,v']) F^{-1}=\ad_{\mk}[Fv,Fv'], \,\forall v,v'\in \mk\},
$$
it is the subgroup of $O(\mk)$ which leaves invariant $\ad_{\mk}([\cdot,\cdot]_{\mid \mk\times\mk})\in(\Lm^2 \mk^*)\otimes\so(\mk)$. Its Lie Algebra 
$$
\der(\mk)=\{A\in \so(\mk)\mid \left[A,\ad_{\mk}[v,v']\right]=\ad_{\mk}[Av,v'] + \ad_{\mk}[v,Av'],\, \forall v,v'\in\mk\}  
$$
is the Lie subalgebra of $\so(\mk)$ which (acting by derivation) leaves invariant $\ad_{\mk}([\cdot,\cdot]_{\mid \mk\times\mk})\in(\Lm^2 \mk^*)\otimes\so(\mk)$.\\
\end{defn}
\begin{thm}\label{m}
Let $M$ be a $G$-symmetric space, $G\subset\Is(M)$, and $\tau\colon\g
\to\g$ an order four automorphism such that $\tau^2=\sigma$.
Then $\tm\in\aut(\mk)$ and $\tau$ can be extended in an unique way to the Lie algebra $\der(\mk)\oplus\mk$ endowed with the Lie bracket
$$
[(A,v),(A',v')]=([A,A']+\ad_{\mk}[v,v'],\, A.v'-A'.v)
$$
and of which $\g$ is a Lie subalgebra, under the inclusion
$a+v\in\hk\oplus\mk\mapsto (\ad_{\mk}a,v)\in\der(\mk)\oplus\mk$, by
setting 
\begin{equation}\label{tau}
\underline{\tau}_{|\mk}=\tm\quad  \text{and}\quad \underline{\tau}_{|\der(\mk)}=\Ad\tm\, .
\end{equation}
Conversely, given $\taum\in O(\mk)$, the linear map $\unt$ defined by (\ref{tau}) is an automorphism of the Lie algebra $\der(\mk)\oplus\mk$ \iif $\taum\in \aut(\mk)$. Besides it 
satisfies $\unt^2 = \Id_{\der(\mk)}\oplus -\Id_{\mk}$  (and in particular is of order four) \iif $\taum\in\Sigma(\mk)$.\\
Hence, to define a locally 4-symmetric bundle over the
Riemannian symmetric space $M$ (with the realisation $M=G/H$, i.e. $\tau$ is an automorphism of $\g$ such that $\tau^2 =  \sigma$) is equivalent to give ourself 
$\tau_{\mk}\in\Sigma(\mk)\cap\aut(\mk)$ such that the order four automorphism $\unt$ of $\der(\mk)\oplus\mk$ stabilizes $\g=\hk\oplus\mk$, i.e. such that 
$\taum(\ad_{\mk}\hk)\taum^{-1}=\ad_{\mk}\hk$ (i.e. $\ad_{\mk}\hk$ is a subalgebra of $\der(\mk)$ invariant by $\Ad\taum$).
Then $\tau=\underline{\tau}_{|\g}$ is an order four automorphism of
$\g$ such that $\tau^2=\Id_{\hk}\oplus -\Id_{\mk}=\sigma$.
\end{thm}
\textbf{Proof.} First $\tm\in\aut(\mk)$: that follows from the fact that $\tau$ is an automorphism, so $\tau\circ\ad a\circ\tau^{-1}=\ad\tau(a)\,, \forall a\in \g$.\\
Second, $\der(\mk)\oplus\mk$ is a Lie subalgebra . We have to check that the Jacobi identity is satisfied. It is a straightforward computation (see \cite{Hel}).
 Then we have to check that $\unt$ is an automorphism \iif $\taum\in\aut(\mk)$.\\
If $\taum\in\aut(\mk)$ then
\begin{description}
\item[$\bullet$]if $A,A'\in\der(\mk)$, $\unt([A,A'])=[\unt(A),\unt(A')]$ because $\unt_{\der(\mk)}=\Ad\taum$ is an automorphism of $\der(\mk)$.
\item[$\bullet$]if $A\in\der(\mk), \, v\in \mk$, $\unt([A,v])=\taum(A.v)=\taum A\taum^{-1}(\taum .v)=[\unt(A),\unt(v)]$
\item[$\bullet$]if $v,v'\in\mk$, $\unt([v,v'])=\Ad\taum(\ad_{\mk}[v,v'])=\ad_{\mk}([\taum v,\taum v'])=[\unt(v),\unt(v')]$ because $\taum\in\aut(\mk)$.
\end{description}
Finally $\unt$ is an automorphism and the unique extension of $\tau$ (because it is determined by $\tm$, see remark~\ref{tau_h}).\\
Conversely if $\unt$ is an automorphism of Lie algebra then
$$
\taum\ad_{\mk}([v,v'])\taum^{-1}=(\unt\,\ad([v,v'])\unt^{-1})_{\mid\mk}=(\ad\unt([v,v']))_{\mid\mk}=\ad_{\mk}([\unt(v),\unt(v')])=\ad_{\mk}([\taum v,\taum v']).
$$
Thus $\taum\in\aut(\mk)$.\\
The last assertion of the theorem follows from what precedes. This completes the proof.\hfill$\blacksquare$

\begin{rmk}\emph{
Let $\taum\in \Sigma(\mk)$ then the condition $\Ad\taum(\ad_{\mk}\hk)=\ad_{\mk}\hk$ implies
that there exists an automorphism $\tauh$ of $\hk$ defined by $\forall a\in \hk$, 
$\Ad\taum(\ad_{\mk}a)=\ad_{\mk}\tauh(a)$, i.e. $\tauh=\ad_{\mk}^{-1}\circ\Ad\taum\circ\ad_{\mk}$.
Then the condition $\taum\in\aut(\mk)$ is equivalent to
$$
\tauh([v,v'])=[\taum v,\taum v'],\ \forall v,v'\in\mk.
$$
And obviously, if these two conditions are satisfied then we have
$\tauh=\tmh$ (where $\tau=\unt_{\g}$ is given by the theorem \ref{m}).
}\end{rmk}
\begin{rmk}\label{s}\emph{
Let us consider the map
$$
s\colon g\in\Is_{p_0}(M)\mapsto\Ad_{\mk}g\circ\tm\circ\Ad_{\mk}g^{-1}\in \Sigma(\mk)
$$
and set $\underline{G_0}=\{g\in\Is_{p_0}(M)| s(g)=\tm\}$. Then $\iso(M)$
acts on $\Sigma(\mk)$ by $g.J=\Ad_{\mk}g\circ J\circ\Ad_{\mk}g^{-1}$
and $s(g)=g.\tm$, and $\underline{G_0}=\mrm{Stab}_{\iso(M)}(\tm)$. In the
same way, the subgroup $H=\iso(M)\cap G$ acts on $\Sigma(M)$ and
$G_0=\mrm{Stab}_H(\tm)$. Then $s(\iso(M))=\iso(M)/\underline{G_0}$ is
a compact submanifold  of $\Sigma(\mk)$, and $s(H)=H/G_0$ is a
relatively compact (immersed) submanifolds of $\Sigma(\mk)$.
}\end{rmk}
\subsubsection{Second convenient hypothesis.}
An other convenient hypothesis on $G$ is to consider that it is a
closed subgroup of $\Is(M)$ (and not only an immersed subgroup).
It is always possible to work with this hypothesis. Let us make
precise this point. Let $\sigma_{p_0}$ be the symmetry of $M$
around $p_0$ (defined by $\sigma$): $\sigma_{p_0}\in\Is(M)$, $\sigma_{p_0}(p_0)=p_0$ and
$T_{p_0}\sigma_{p_0}=-\mathrm{Id}$. Then $\sigma_{p_0}$ belongs to
the isotropy subgroup $\Is_{p_0}(M)=\{f\in \Is(M); f(p_0)=p_0\}$,
and we can define the involution of $\Is(M)$:
$$
\sigma_{\Is(M)}=\mathrm{Int}(\sigma_{p_0})\colon g\in\Is(M)\mapsto
\sigma_{p_0}\circ g\circ\sigma_{p_0}^{-1}\in \Is(M).
$$
It is easy to see that we have
\begin{equation}\label{fix}
(\Is(M)^{\sigma_{\Is(M)}})^0\subset
\Is_{p_0}(M)\subset \Is(M)^{\sigma_{\Is(M)}}
\end{equation}
(see \cite{Hel,besse}). The result of this is that $\sigma\colon G\to G$ is the
restriction of $\sigma_{\Is(M)}$ to $G\subset \Is(M)$ (they induce $\sigma_{p_0}$
on $M=G/H$ and the identity on $H$, thus, since $G$ is locally isomorphic to
$M\times H$, they are identical, see also \cite{Hel}). Moreover there exists
an unique subgroup $\Bar{G}$ of $\mathrm{Diff}(M)$ such that for any
$G$-invariant Riemannian metric $b$ on  $M$, the group $\bar{G}$ is
the closure of $G$ in $\Is(M,b)$: $\Is(M,b)$ is closed in $\mrm{Diff}(M)$ and so the
closure of $G$ in $\Is(M,b)$ is its closure in $\mrm{Diff}(M)$ and thus it does not depend on $b$
(see \cite{besse,Hel}). Then $\sigma$  extends in an unique way to an
involutive morphism $\bar{\sigma}\colon\Bar{G}\to\Bar{G}$, which is the restriction of
$\sigma_{\Is(M)}$ to $\Bar{G}$. Hence denoting by $\hat{H}$
the isotropy subgroup of $\Bar{G}$ at $p_0$, $\hat{H}=\Is_{p_0}(M)\cap\Bar{G}$,
we have according to (\ref{fix}), $(\Bar{G}^{\bar{\sigma}})^0\subset\hat{H}
\subset\Bar{G}^{\bar{\sigma}}$. Besides $\bar{\sigma}$ gives rise to the symmetric
decomposition $\lie \Bar{G}=\lie\hat H\oplus \mk$.\\
In addition to that, we have $\hat{H}=\Bar{H}$. Indeed, let
$\Phi\colon U\times\Is_{p_0}(M)\to \Is(M)$ be a local trivialisation
of $\Is(M)\to M$, such that $\Phi(p_0,h)=h$, and $\Phi(U\times H)=\Phi(U\times
\Is_{p_0}(M))\cap G$ (take $\Phi(p,h)=\phi(p).h$, with $\phi\colon U\to G$ a local
section such that $\phi(p_0)=1$). Further, if $ g\in \Is_{p_0}(M)\cap \Bar{G}$ and
$(g_n)$ is a sequence of $G\cap \Phi(U\times \Is_{p_0}(M))$ such that $g_n\to g$,
then $\Phi^{-1}(g_n)=(u_n,h_n)\in U\times H$ converges to
$\Phi^{-1}(g)=(p_0,g)$, thus $h_n\to g$ so $g\in \Bar H $.\\
Moreover, $\Bar H$ is a closed subgroup of $\iso(M)$, thus it is
compact. Hence, we have the symmetric realisation $M=\Bar G/\Bar H$ and
$\Ad_{\mk}(\Bar H)$ is compact: we have showed that the hypothesis
"$\Ad_{\mk}(H)$ relatively compact" and "$\Ad_{\mk}(H)$ compact" give the same symmetric spaces. Moreover,  by using the preceding reasoning (to prove $\hat{H}=
\Bar H$) it is easy  to see that if $\Ad_{\mk}(H)$ is compact then $G$ is
closed in $\Is(M)$ (see also \cite{Hel}) so that the hypothesis "$\Ad_{\mk}(H)$ is compact" and "$G$ is
closed in $\Is(M)$" are in fact equivalent.\\
Besides, the closure of
$G$ is  the same in $\Is(M)$ and in $\Is(\tilde{M})$ with
$\tilde{M}=G/H^0$: since $M$ and $\tilde{M}$ are complete (a
Riemannian homogeneous space is complete) then $\Is(M)$ and
$\Is(\tilde{M})$ are complete (see \cite{Hel}), and thus the closure of $G$ in one
of this group is the completed of $G$. \\
Now, let us suppose that we have a locally 4-symmetric bundle over
$M$.
\begin{thm}\label{bar}
Let us consider $M$ a $G$-symmetric space with $G\subset\Is(M)$ and
$\tau\colon\g\to\g$ an order four automorphism such that $\tau^2=\sigma$.
Then the extension $\unt$ of $\tau$, given by theorem \ref{m} stabilizes the Lie algebra, $\lie\Bar G$, of the closure of $G$ in $\Is(M)$:
$$
\unt(\lie\Bar G)=\lie \Bar G.
$$
Then denoting by $\bar\tau :=\unt_{\mid\lie\Bar G}$ the extension of $\tau$ to $\lie\Bar G$ (given by 
theorem~\ref{m}), the subgroup fixed by $\bar\tau$ (defined by (\ref{g-0})) is the closure of $G_0$:
$$
\Hat{G}_0:=\{g\in \Bar H\mid\tau_{|\mk}\circ\Ad_{\mk} g\circ\tau_{|\mk}^{-1}=\Ad_{\mk} g\}=\Bar{G}_0.
$$
Finally the new locally 4-symmetric bundle over $M$ defined by $\bar\tau$ is $\Bar {G}/\Bar{G_0}$, and using the notation of
remark~\ref{s}, the fibre of $\Bar {G}/\Bar{G_0}$, $\Hat{S}_0:=s(\Bar H)=\Bar H/\bar{G}_0$, is the closure of the fibre $S_0$ of $G/G_0$, $\overline{S_0}=\overline{s(H)}=\overline{(H/G_0)}$, in the maximal fibre
over $M$: $\und{S_0}:=s(\iso(M))=\iso(M)/\und{G_0}$.
\end{thm}
\begin{rmk}\emph{
In particular, if we suppose that we have an order four automorphism $\tau$
of $G$, such that $\tau^2=\sigma$, then since $\tau$ is uniformly continuous,
it extends into an order four automorphism $\bar{\tau}\colon\Bar{G}\to\Bar{G}$ (because
$\Is(M)$ is complete) and obviously $\bar{\tau}^2=\bar{\sigma}$.}
\end{rmk}
The following theorem precises the link between the Lie algebra setting and the one of the Riemannian symmetric space $M$ (first point of theorem~\ref{aut*}), which will allow us (in theorem~\ref{J_0}) to translate the theorem~\ref{m} in terms of the setting of $M$. The two last points (of theorem~\ref{aut*}) characterize the "satisfying cases": any element in $\aut(\mk)$ defines an automorphism in $\lis(M)$ (an example of "unsatisfactory" case is given by $M=\R^{2n-r}\times \mathbb{T}^r$, see section~\ref{euclidean}).
\begin{thm}\label{aut*}
Let us consider $M$ a  Riemannian symmetric space and $\tl M$ its universal covering.\\ 
$\bullet$ The curvature operator (in $M$) is given by $R_{p_0}(\cdot,\cdot)=-\ad_{\mk}([\cdot,\cdot]_{\mid \mk\times\mk})$ and thus\footnote{$\mrm{Hol}(M)$ is the holonomy group of $M$}
\begin{equation}\label{is1}
\begin{array}{l}
\der(\mk)=\liso(\tl M)\supset\liso(M)\supset\lie\mrm{Hol}(M)\\
\aut(\mk)\supset\iso(\tl M)\supset\iso(M)\supset\mrm{Hol}(M)
\end{array}
\end{equation}
(using the identification $T_{p_0}M=\mk$) and $\der(\mk)\oplus\mk=\lis(\tl M)$.\\
$\bullet$ Moreover the following statements are equivalent:\\
(i) $\liso(\tl M)=\liso(M)$ (i.e. $\lis(\tl M)=\lis(M)$)\\
(ii) $M=M'\times M_0$, with $M'$ of the semisimple type (i.e. $\Is(M')$ is semisimple) and $M_0$ is Euclidean.\\
(iii) $\hk_0=\so(\mk_0)$, where $\hk_0$ and $\mk_0$ are respectively the Euclidean part of $\liso(M)$ and $\mk$ respectively, in the decomposition $\lis(M)=\g'\oplus \g_0$, with $\g'$ semisimple and $\g_0$ of the Euclidean type.\\[1mm]
$\bullet$ Moreover the following statements are also equivalent:\\
(i) $\liso(\tl M)=\liso(M)\oplus\so(\mk_0)$\\
(ii) $\hk_0=0$\\
(iii) Let $\tl M= \tl M'\times \tl M_0$ be the decomposition of $\tl M$ into the semisimple  and Euclidean type, $\Gamma$ the group of deck transformations of the covering $\pi\colon\tl M\to M$. Then  the projection on the Euclidean factor (in $\Is(\tl M)=\Is(\tl M')\times\Is(\tl M_0)$) of $\Gamma$ satisfies $\Gamma_0\cong\Z^r$ with $r=\dim \tl M_0$ so that $\tl M_0/\Gamma_0=\mathbb{T}^r$.\\[1mm]
Further $\aut(\mk)$ stabilizes $\liso(M)$ \iif one of the 6 previous statements holds i.e. \iif $\lis(\tl M)/\lis(M)=\{0\}$ or $\so(\mk_0)$. Denoting by $\aut(\mk)^*$ the subgroup of $\aut(\mk)$ which stabilizes $\liso(M)$,  then  the maximal subalgebra of $\liso(M)$ invariant by $\aut(\mk)$ is $\liso(M)$  if $\aut(\mk)=\aut(\mk)^*$ and $\hk'=\liso(\tl M')$ if not.
\end{thm}

Let us consider $M$ a $G$-symmetric space with $G\subset\Is(M)$ and
$\tau\colon\g\to\g$ an order four automorphism such that $\tau^2=\sigma$.\\ 
 Then the extension $\unt$ of $\tau$, given by theorem \ref{m} defines a maximal locally 4-symmetric bundle over $M$.
Indeed let $\und\g$ be the maximal subalgebra of $\lis(M)$ invariant by $\unt$, and $\und G$ the subgroup of $\Is(M)$ generated by it. Then  $\und G$ is a closed subgroup of 
$\Is(M)$ acting symmetrically on $M$: since  $\und G$ is connected, it is invariant by $\sigma_{\Is(M)}$ (and it contains $G$) thus it acts  symmetrically on $M$, then $\und G$ is closed as an immediate consequence of the maximality and theorem~\ref{bar}. Therefore, $\unt_{\mid\und \g}$ defines a maximal locally 4-symmetric bundle over $M$, with the realisation $M=\und G/\und H$.\\ 
We can also define a minimal locally 4-symmetric bundle over $M$, by considering the subalgebra $\g'\oplus\mk_0$ (where $\g'$ is the semisimple part of $\lis(M)$ and  $\mk_0$ the 
Euclidean part of $\mk$).
\begin{thm}\label{J_0}
In conclusion, given any (even-dimensional) Riemannian symmetric space $M$, to define over it a locally 4-symmetric
bundle is equivalent to give ourself $J_0\in\Sigma(T_{p_0}M)$ which leaves invariant the curvature $R_{p_0}$:
$$
R_{p_0}(J_0 X,J_0 Y)=J_0R_{p_0}(X,Y)J_0^{-1}\ \forall X,Y\in T_{p_0}M.
$$
 Then the order four automorphism of $\mathfrak{Is}(\tl M)$, $\unt$, 
defined by $J_0$, defines the maximal locally 4-symmetric bundle over $M$, $\und{N}=\und G/\und{G_0}$ with 
$\und{G_0}=\{g\in\und H\mid J_0 g J_0^{-1}= g\}$. Moreover, any locally 4-symmetric bundle over $M$ is a subbundle of one such maximal 
bundle and to define such a  subbundle $N$ is equivalent to give ourself a Lie subgroup $G\subset \Is(M)$ acting symmetrically on $M$ such that $\unt(\g)=\g$ i.e. $\Ad J_0(\hk)=\hk$.
In this case, the closure $\Bar{N}$ of $N=G/G_0$ in the (unique) maximal locally 4-symmetric bundle over $M$ containing $N$, $\und{N}$, is also a locally 4-symmetric bundle over 
$M$ and we have $\Bar N=\Bar G/\Bar{G_0}$, $M=\Bar G/\Bar H$.
\end{thm}
\textbf{Proof of theorem~\ref{bar}}  
We have to check that $\underline{\tau}(\lie\Bar G)=\lie\Bar G$, i.e.
according to the theorem \ref{m}, $\Int\taum(\lie\Bar H)=\lie\Bar
H$. We still have $\Int\taum(H^0)=H^0$, thus $\Int\taum(\overline{H^0})=\overline{H^0}$. It remains to
verify that $(\Bar H)^0=\overline{H^0}$. But this is simply the result of the
fact that $\tilde{M}:=G/H^0=\Bar G/\overline{H^0}$ (the closures in $\Is(\tl M)$ and $\Is(M)$ are the same) is a discrete covering of $M=\Bar G/\Bar H$. Indeed $(\Bar H)^0$ is closed thus $(\Bar H)^0\supset \overline{H^0}$
and then we have
$$
\tl{M}=\Bar{G}/\overline{H^0}\xrightarrow{\text{fibration}}\Bar{G}/(\Bar H)^0
\xrightarrow{\text{covering}}\Bar G/\Bar H =M
$$ 
 and $\tl M \xrightarrow{\text{covering}} M$, hence $(\Bar H)^0/\overline{H^0}$
is discrete but the two groups are connected ($(\Bar H)^0$ suffices) thus
$(\Bar H)^0=\overline{H^0}$. We have proved that $\unt(\lie\Bar G)=\lie \Bar G$.\\
Using the notation of remark \ref{s}, we have, since $\Bar H$ is
compact, $s(\Bar H)=\overline{s(H)}$, hence using the same method as for
$\Hat{H}$, we can easily show that
$\Hat{G_0}:=s^{-1}(\tm)\cap\Bar H=\Bar G_0$ and thus $s(\Bar
H)=\Bar H/\Bar G_0$. Finally, the new locally 4-symmetric space is
$\Bar G/\Bar G_0$. This completes the proof.\hfill$\blacksquare$\\
\textbf{Proof of theorem~\ref{aut*}} For the first point see \cite{Hel}. For the following points, see sections~\ref{semisimple} and \ref{euclidean}.\hfill$\blacksquare$\\
\textbf{Proof of theorem~\ref{J_0}} The first assertions are nothing but the translation of theorem~\ref{m}, using theorem~\ref{aut*}. 
Then, we have to prove that $\Bar G/\Bar G_0$ is the closure $\Bar N$ of $N=G/G_0$ in $\und G/\und{G_0}$.
Let $\pi_{J_0}\colon\und G\to\und G/\und{G_0}$ be the projection map, then we have $\pi_{J_0}(G)=G/\und{G_0}\cap G
=G/G_0$ (according to definition (\ref{g-0})) and thus $\pi_{J_0}(\Bar G)\subset\overline{\pi_{J_0}(G)}=\Bar N$ but 
 $\pi_{J_0}(\Bar G)=\Bar G/\und{G_0}\cap \Bar G =\Bar G/\Bar{G}_0$ (according to definition (\ref{g-0}) and $\Hat{G}_0=\Bar G_0$). Hence $\Bar G/\Bar G_0\subset \Bar N$. These are together subbundle (over $M$) of $\und N$ and using a trivialisation of $\und N=\und G/\Bar G_0\to M$ (same reasoning as for $\Hat H$) it is easy to see that  the fibre of $\Bar N$ (over $p_0$) is $\Bar H/\Bar G_0$ which implies that $\Bar G/\Bar G_0 = \Bar N$. This completes the proof.\hfill$\blacksquare$
\begin{rmk}\emph{
According to the definition (\ref{g-0}), $\tm$ and $-\tm$ give rise to the same group $G_0$. Moreover $\tm=(\tau^{-1})_{\mid\mk}$ and in particular if $\tau$ integrates in $G$ then $G^{\tau}=G^{\tau^{-1}}$. Besides $(\tau^{-1})^2=\sigma^{-1}=\sigma$, hence $\tau^{-1}$ defines the same locally 4-symmetric bundle over $M$ as $\tau$.
Moreover, given any $\taum\in\Sigma(\mk)\cap\aut(\mk)$, then $-\taum\in\Sigma(\mk)\cap\aut(\mk)$
and gives rise (according to theorem \ref{m}) to the automorphism $\unt^{-1}$ which gives rise to the same maximal locally 4-symmetric bundle over $M$ and the same family of 4-symmetric subbundle over $M$.
}\end{rmk}
From now, we will always suppose that $G$ is a closed subgroup of
$\Is(M)^0$. The result of this is that the isotropy subgroup of $G$ at the point
$p_0$, $H=\mathrm{Stab}_G(p_0)$ is  compact  and can be identified (via the adjoint
representation on $\mk$, resp. via the linear isotropy representation)  to a closed subgroup of $O(\mk)$ (resp. of $O(T_{p_0}M)$). Then according to theorem \ref{J_0}, to study the case of non-closed subgroup of $\Is(M)^0$ (or equivalently the non-closed locally 4-symmetric bundle over $M$), we have just to consider the non-closed subgroups of our closed group $G$, acting symmetrically on $M$, and whose Lie algebra is invariant by $\tau$.
%
%%%%%%%%%%%%%%%%%%%%%%%%%%%%%%%%%%%%%%%%%%%%%%%%%%%%%%%%%%%%%%%%%%%%%%%%%%%%%%%%%%%%
%
\subsection{Twistor subbundle}\label{twistor}
We give ourself a locally 4-symmetric bundle $N=G/G_0$ (defined by an order four automorphism $\tau$ and by (\ref{g-0})) over a symmetric
space $M=G/H$. We will show that $G/G_0$ is a subbundle of the twistor bundle
$\Sigma(G/H)$ .
Under the isomorphism between $TM$ and $[\mk]=\{(g.p_0,\Ad g(\xi)),\, \xi\in\mk ,g\in
G\}$,
$T_{p_0}M$ is identified to $\mk$: $\xi\in\mk\mapsto\xi .p_0\in T_{p_0}M$ is
an isomorphism of vector spaces.
Then to $\tau_{|\mk}\in\Sigma(\mk)$ (resp. to $-\tm=\tm^{-1}\in\Sigma(\mk)$) corresponds $J_0\in\Sigma(T_{p_0}M)$, and more
generally to $\Ad g\circ\tau_{|\mk}\circ\Ad g^{-1}\in\Sigma(\Ad g(\mk))$ (resp. $\Ad g\circ\tau_{|\mk}^{-1}\circ\Ad g^{-1}$)
corresponds $gJ_0g^{-1}\in\Sigma(T_{g.p_0}M)$. Thus we have
defined a map
$$
  \begin{array}{crcl}
\rho_{J_0}\colon & G & \longrightarrow &\Sigma(M)\\
 & g & \longmapsto & gJ_0g^{-1}
  \end{array}
$$
which according to the definition (\ref{g-0}) of $G_0$ gives rise under quotient to
the injective map:
$$
 \begin{array}{crcl}
i \colon & G/G_0 & \longrightarrow &\Sigma(M)\\
 & g.G_0 & \longmapsto & gJ_0g^{-1}\ .
  \end{array}
$$
Moreover, $i$ is an embedding. Indeed, $G$ acts smoothly on the
manifolds $\Sigma(M)$ and so the map $g\in G\mapsto gJ_0g^{-1}\in
\Sigma(M)$ is  of constant rank. Thus
$i\colon G/G_0\to\Sigma(M)$ is an injective map of constant rank and so it
is an injective immersion. We can add that $i\colon
G/\mathrm{Stab}_G(J_0)\to G.J_0$ is an homeomorphism if the orbit
$G.J_0$ is locally closed in $\Sigma(M)$ (see \cite{dieudonne}).  
We will show
directly that $i(G/G_0)=G.J_0$ is a subbundle of $\Sigma(M)$.\\
First, let us precise the fibration $G/G_0\to G/H$ . We have the
isomorphism of bundle: $G/G_0\cong G\times_{H}H/G_0$. In particular,
the fibre type of $G/G_0$ is $H/G_0$. Besides $i$ is a morphism of bundle (over $M$).
 Since $i$ is also an
injective immersion, we can identify the fibres of $G/G_0$ and
$i(G/G_0)$ respectively over the point $g.p_0\in M$. The fibre of $i(G/G_0)$ over $p=g.p_0$ is
$gS_0g^{-1}$ where $S_0=\mathrm{Int}(H)(J_0)\subset\Sigma(T_{p_0}M)$ is the fibre
over $p_0$.\footnote{we remark that $H\subset O(T_{p_0}M)$, $G_0\subset U(T_{p_0}M,J_0)$ 
and $S_0=H/G_0$ is a compact submanifold of $\Sigma(T_{p_0}M)$.}\\
Now let us show that $i(G/G_0)$ is a subbundle of $\Sigma(M)$. Let
$\sigma\colon U\subset G/H\to G$ be a local section of the fibration
$\pi_H\colon G\to G/H$. Then we have the following trivialisation of
$\Sigma(U)$:
$$
\Phi\colon (p,J)\in
U\times\Sigma(T_{p_0}M)\longmapsto (p,\sigma(p)J\sigma(p)^{-1})\in\Sigma(U)
$$
and we have $\Phi(U\times S_0)=\bigsqcup\limits_{p\in U}\{p\}\times(\sigma(p)S_0\sigma
(p)^{-1})=i(G/G_0)\cap\Sigma(U)$. Thus  $i(G/G_0)$ is a subbundle
of $\Sigma(M)$, hence $i$ is an embedding.

Let us recapitulate what precedes:

\begin{thm}
The map
$$
 \begin{array}{crcl}
i \colon & G/G_0 & \longrightarrow &\Sigma(M)\\
 & g.G_0 & \longmapsto & gJ_0g^{-1}\ .
  \end{array}
$$
is an embedding and a morphism of bundle from $G/G_0$ into
$\Sigma(M)$. Besides the fibre of $i(G/G_0)$ over the point
$p=g.p_0$ is  $gS_0g^{-1}$, with $S_0=\mathrm{Int}(H)(J_0)$ and
$J_0\in\Sigma(T_{p_0}M)$ corresponding to $\tau_{|\mk}\in\Sigma(\mk)$ (resp.
to $\tau_{|\mk}^{-1}$) .
\end{thm}
\begin{rmk}\label{explanation}\emph{
Given one order four automorphism, we have two different ways to embed $G/G_0$ into $\Sigma(M)$ by $J_0=\pm\tm$. The two submanifolds $i_{J_0}(G/G_0)$ and $i_{-J_0}(G/G_0)$ are isomorphic by $J\mapsto -J$. These are identical \iif $H$ contains  an element which anticommutes with $J_0$. If $\dim M=2$ mod $4$ then they lie in  different connected components of the twistor space (one in $\Sigma^+(M)$ and the other one in $\Sigma^-(M)$). In theorem~\ref{fondamental} we use $-\tm$ to respect the convention: $\alpha_{-1}$ is an $(1,0)$-form.}
\end{rmk}
\begin{rmk}\emph{
If  we consider a
locally 4-symmetric bundle $N=G/G_0$ over $M$, with $G_0$ not defined by (\ref{g-0}), then $i$ is not injective in general: to obtain
an injective map $i$, we must consider the locally 4-symmetric space
$G/\pi_K^{-1}(G_0')=(G/K)/G_0'$ where $K=\ker \Ad_{\mk}$ and $G_0'$
is the subgroup of $H'=H/K$ defined by (\ref{g-0}). In particular, we see that in
general a 4-symmetric space $G/G_0$ is not a submanifold of a twistor space
(see section~\ref{example}). Moreover, we can
see the aim of our definition (\ref{g-0}) (and in particular of our convention $G_0=G^{\tau}\cap H$): it is to obtain an injective map $i$.}
\end{rmk}
\begin{rmk}\label{tl-g}\emph{
For any covering $\pi\colon\tl G\to G$, $\tl G$ acting symmetrically on $M$, we have $\tl\rho_{J_0}(\tl G)=\rho_{J_0}(G)=i_{J_0}(G/G_0)$ : the locally 4-symmetric subbundle of $\Sigma(M)$, $i_{J_0}(G/G_0)$ does not depend on the choice of the group $G$ (we have chosen for $G$, the subgroup of $\Is(M)$ generated by $\g$).\\
\indent Moreover, $\rho_{J_0}(\und G)=i_{J_0}(\und N)$ is a maximal locally 4-symmetric subbundle in $\Sigma(M)$. Now, suppose that we work with a non-closed subgroup $G'\subset \Is(M)$, then $\rho_{J_0}(G')=i_{J_0}(G'/G_0')$ is an immersed subbundle in $\Sigma(M)$: $\Phi(U\times S_0')=\bigsqcup\limits_{p\in U}\{p\}\times(\sigma(p)S_0'\sigma
(p)^{-1})=i(G'/G_0')\cap\Sigma(U)$ but the fibre $S_p'=\sigma(p)S_0'\sigma(p)^{-1}$ is only a (non-closed relatively compact) immersed submanifold in $\Sigma(T_p M)$. And since $i_{J_0}$ is an embedding (from $\und{N}$ into 
$\Sigma(M)$) we have $i_{J_0}(\Bar G'/\Bar G_0')=i_{J_0}(\Bar N')=\overline{i_{J_0}(N')}$. In others words,
taking the closure of $G'$ in $\Is(M)$ is equivalent to take the closure of $N'$ in $\und N$ according to theorem~\ref{J_0} which is equivalent to take the closure of $i_{J_0}(N')$ in $i_{J_0}(\und N)$.
}\end{rmk}
\begin{rmk}\label{tlg}\emph{
The maximal locally 4-symmetric bundles $\und N$ are disjoint : these are orbits in $\Sigma(M)$. More precisely these are suborbits of $\Is(M)^0$-orbits in the form $\und G\cdot
J_0$ in $\Sigma(M)\cap \aut(M)$ with $\aut(M)=\bigsqcup_{p\in M} \{p\}\times \aut(T_p M)$ (see Sections~\ref{semisimple} and \ref{euclidean}). In particular, $\Sigma(M)\cap \aut(M)$
is the disjoint union of all the maximal locally 4-symmetric bundles over $M$. Moreover, the set of maximal locally 4-symmetric bundles over $M$, $\mathcal N(M)$, contains the subset
$\mathcal N^*(M)$ of elements which are $\Is(M)^0$-orbits, i.e. elements $\rho_{J_0}(\Is(M)^0)$  with $J_0\in\Sigma(\mk)\cap \aut^*(\mk)$ :  $\mathcal N^*(M)=(\Sigma(M)\cap \aut^*(M))/\Is(M)^0\subset  \Sigma(M)/\Is(M)^0$.
}\end{rmk}
\begin{rmk}\label{sigmatl}\emph{
The Riemannian manifold $M=G/H$ is orientable \iif
$\Ad_{\mk}H\subset SO(\mk)$ (or equivalently $H\subset SO(T_{p_0}M)$).
Besides $\tau_{|\mk}\in\Sigma^{\eps}(\mk)$, and to fix ideas, let us
suppose that $\eps=1$. Then, if  $M$ is orientable, $i$ is an
embedding from $G/G_0$ into $\Sigma^+(M)$.
Moreover, if we work with $\tilde{M}=G/H^0$, we are sure that
$H^0\subset SO(T_{p_0}\tilde{M})$. Hence, if we work with what we
called the first possibility (see section~\ref{first}), $i$ takes values in $\Sigma^+(\tl M)$. In
other words, given a locally 4-symmetric bundle over $M$, the
corresponding 4-symmetric bundle over its universal covering $\tl M$
(see remark \ref{univcov}) is embedded in $\Sigma^+(\tl M)$.\\
Let us consider more generally any covering $\pi\colon\tilde{M}\to M$ then
it induces the covering
$\pi_{\Sigma}\colon\Sigma(\tl{M})\to\Sigma(M)$ which is also a
morphism of bundle over $\pi\colon\tl{M}\to M$. It is defined by
$$
\pi_{\Sigma}\colon J_{\tl x}\in\Sigma(T_{\tl x}\tl M)\mapsto T_{\tl
x}\pi\circ J_{\tl x}\circ (T_{\tl x}\pi)^{-1}\in \Sigma(T_{x}M).
$$
Now, let us suppose that $\pi$ comes from a covering $\tl \pi\colon
\tl G\to G$ and that we have $M=G/H$, $\tl M=\tl G/\tl{H}^0$
(symmetric realisation) with $\tl H=\tl\pi^{-1}(H)$ and $G\subset\Is(M)$, $\tl G\subset\Is(\tl M)$ 
(see section~\ref{first}). 
Then we have
$$
T_{\tl x}\pi\circ(\tl g J_{\tl p_0}\tl g^{-1})\circ (T_{\tl x}\pi)^{-1}=
g J_{p_0} g^{-1}
$$
with $\tl x=\tl g.\tl p_0$, $g=\tl\pi(\tl g)$. Hence the restriction
of $\pi_{\Sigma}$ to $\tl G/\tl G_0$ gives rise to the morphism of
bundle (\ref{diag})\footnote{i.e. $\pi_{\Sigma}\circ i_{J_{\tl p_0}}=
i_{J_{p_0}}\circ\pi_0$, where $\pi_0\colon \tl G/\tl G_0\to G/G_0$ is given by (\ref{diag}).}. Moreover\footnote{In fact, $h\tl S_0 h^{-1}$ means obviously
$T_{h.p_0}\pi\circ(h\tl S_0 h^{-1})\circ T_{h.p_0}\pi^{-1}$. $\pi_{\Sigma}$ allows to consider the fibres
$\Sigma(T_{\tl x}\tl M)$ as included in the fibre $\Sigma(T_{x}M)$, with $x=\pi(\tl x)$.}
$$
S_0=\Int(H)(J_0)=\bigcup_{h\in\tl H/\tl{H}^0}h\tl S_0 h^{-1}
$$
with, since $\tl{H}^0\subset SO(T_{p_0}\tl M)$, $\tl S_0\subset\Sigma^+(T_{p_0}\tl M)$.
Further if $H\subset O(T_{p_0}M)$ is not included in $SO(T_{p_0}M)$ (i.e. $M$ is not
orientable), then we have
$$
\pi_{\Sigma}(\Sigma^+(\tl M))=\Sigma(M).
$$
}\end{rmk}
\begin{rmk}\emph{
Let us see what happens when we change $M$, in theorem~\ref{fondamental}. Let $\tl G$ be a covering of $G$, acting symmetrically on a covering $\tl M$ of $M$, $\pi\colon\tl M\to M$, with $\tl G\subset \Is(\tl M)$. Then according to 
remark~\ref{sigmatl}, we have $\pi_{\Sigma}\circ i_{J_{\tl p_0}}=i_{J_{p_0}}\circ\pi_0$, with $\pi_0\colon \tl G/\tl G_0\to G/G_0$  the morphism of bundle (over $\pi\colon\tl M\to M$) given by (\ref{diag}). Then given any solution $\alpha$ of the $(\g,\tau)$-system~(\ref{courbnul}), let us integrate it in $\tl G$ and $G$ respectively, $\tl U\colon L\to \tl G$, $U\colon L\to G$ with $\tl U(0)=1$, $U(0)=1$ ($0$ is a reference point in $L$), we have $\tl \pi\circ\tl U=U$. Then let us project these lifts in $\tl G/\tl G_0$ and $G/G_0$ respectively: we obtain the geometric solutions  $\tl J\colon L\to \tl G/\tl G_0$ and $J\colon L\to G/G_0$ respectively and we have $\pi_0\circ\tl J=J$. Then let us embed these into the twistor spaces $\Sigma(\tl M)$ and $\Sigma(M)$ to obtain the admissible twistor lifts $\tl J_{\tl X}\colon L\to i_{\tl J_0}(\tl G/\tl G_0)$ and $J_{X}\colon L\to i_{J_0}(G/G_0)$ respectively which are related by $\pi_\Sigma\circ \tl J_{\tl X}=J_{X}$, and in particular
$\pi\circ\tl X=X$.
}\end{rmk}
%
%
%
%%%%%%%%%%%%%%%%%%%%%%%%%%%%%%%%%%%%%%%%%%%%%%%%%%%%%%%%%%%%%%%%%%%%%%%%%%%%%%%%%%%
%
%
%
%%%%%%%%%%%%%%%%%%%%%%%%%%%%%%%%%%%%%%%%%%%%%%%%%%%%%%%%%%%%%%%%%%%                                                                %  
%  SECTION 3: SPLITTING OF M INTO THE 3 TYPE OF SYMMETRIC SPACE  % %%%%%%%%%%%%%
%                                                                %
%%%%%%%%%%%%%%%%%%%%%%%%%%%%%%%%%%%%%%%%%%%%%%%%%%%%%%%%%%%%%%%%%%

\section{Splitting of $M$ into the 3 types of symmetric spaces}\label{splitting}
In the following theorems and corollaries, we study the behaviour of the automorphism $\tau$ with respect to the de Rham decomposition of $M$.
\begin{thm}\emph{\cite{Hel,besse}}
Let $M$ be a simply connected Riemannian symmetric space. Then $M$ is a product
$$
M=M_0\times M_-\times M_+
$$
where $M_0$ is an Euclidean space, $M_-$ and $M_+$ are Riemannian symmetric spaces of the compact and non-compact types respectively. In particular
$$
M=M_0\times M'
$$
where $M'$ has a group of isometries $G=\Is(M')$ semisimple and its isotropy subgroup at $p_0\in M'$, $H$, (which is connected because $M'$ is simply connected) is equal to the 
holonomy group of $M'$. Hence a Riemannian symmetric space $M$ of which the isometry group is semisimple (which is equivalent to say that its universal covering has not Euclidean 
factor, or equivalently the Lie algebra of $G$ does not contain non-trivial abelian ideal, i.e. its Killing form is non-degenerated) has a unique symmetric realisation $G/H$, with $G$ 
acting effectively. In this unique realisation, we have necessarily $G=\Is(M)^0$ \footnote{as usual, we suppose that $G$ is connected} and 
$H=\iso^0(M):=\iso(M)\cap\Is(M)^0(\supset \iso(M)^0)$.
Further the Lie algebra $\mathfrak{Is}_{p_0}(M)=\der(\mk)=\mathfrak{Hol}(M)$ is spanned by $[\mk,\mk]=\{R_{p_0}(X,Y),\, X,Y\in T_{p_0}(M)\}$.\\
Moreover the universal covering of such a Riemannian symmetric space $M$, admits a decomposition into a product of irreducible Riemannian symmetric spaces (i.e with  linear isotropy representations  which are irreducible)
$$
\tl M = M_1\times\cdots\times M_r .
$$
\end{thm}
\begin{thm}\label{alpha}
Let us consider the decomposition of $(\g,\sigma)$ into the sum of orthogonal (for the Killing form) ideals of the compact, non-compact and Euclidean types respectively:
$$
\g=\mak{l}_0\oplus\mak{l}_-\oplus\mak{l}_+
$$
and let $\mak l_{\alpha}=\hk_{\alpha}\oplus\mk_{\alpha}$ be the eigenspace decomposition of the involution
$\sigma_{\mid \mak{l}_\alpha}$.\\
Suppose now that we have an order four automorphism $\tau\colon\g\to\g$ with $\tau^2=\sigma$. Then $\tau(\mak l_{\alpha})=\mak l_{\alpha}$, $\tau(\hk_{\alpha})=\hk_{\alpha}$, 
$\tau(\mk_{\alpha})=\mk_{\alpha}$ for $\alpha=0,-,+$. Hence $\taum=\oplus_\alpha \tau_{\mk_{\alpha}}$, with $\tau_{\mk_{\alpha}}\in\Sigma(\mk_\alpha)$, and 
$\tau_{\mid\mak l_{\alpha}}$ is the automorphism of $\mak l_{\alpha}$ defined by $\tau_{\mk_{\alpha}}$ according to theorem \ref{m} and we have $\tau_{\mid\mak l_{\alpha}}^2=
\sigma_{\mid \mak{l}_\alpha}$. Moreover, we have $\aut(\mk)=\prod_{\alpha}\aut(\mk_\alpha)$.
\end{thm}
\begin{cory}\label{alpha'}
Let $M$ be a $G$-symmetric space, $G\subset\Is(M)$ and $\tau\colon\g\to\g$ an order four automorphism with $\tau^2=\sigma$. Let $\tl M$ be its universal covering, which has a symmetric realisation $\tl M=\tl G/\tl H^0$, with $\pi\colon\tl G\to G$ a covering of $G$, $\tl H=\pi^{-1}(H)$ and $\tl G\subset\Is(\tl M)$, such that $\tau$ integrates into $\tl\tau\colon \tl G\to\tl G$. Then the decomposition of $\g$ into 3 ideals of different types gives rise to the following decomposition of $\tl G$:
$$
\tl G=L_0\times L_- \times L_+
$$
 which induces the following decompositions of $\tl H^0$ and $\tl G_0=\tl H^0\cap\tl G^{\tl\tau}$, corresponding also to the decompositions 
$\hk=\oplus_\alpha \hk_\alpha$ and $\g_0=\oplus_\alpha (\g_0)_\alpha$:
\begin{eqnarray}
\tl H^0 & = & H_0\times H_- \times H_+\\
\tl G_0 & = & (G_0)_0\times (G_0)_- \times (G_0)_+\, .
\end{eqnarray}
Hence $\tl M= M_0\times M_-\times M_+$ and $\tl N= N_0\times N_-\times N_+$ with $M_\alpha=L_\alpha/H_\alpha$,
$N_\alpha=L_\alpha/(G_0)_\alpha$. Besides $\tl\sigma$ and $\tl\tau$ admit the decompositions $\tl\sigma=\prod_\alpha\tl\sigma_\alpha$ and $\tl\tau=\prod_\alpha\tl\tau_\alpha$, and $H_\alpha= (L_\alpha^{\tl\sigma_\alpha})^0$, $(G_0)_\alpha=H_\alpha
\cap L_\alpha^{\tau_\alpha}=(L_\alpha)_0$. Moreover $N_\alpha$ is a 4- symmetric bundle over $M_\alpha$.
\end{cory}
\begin{thm}\label{gi}
Let us consider the decomposition of $(\g,\sigma)$ into the sum of orthogonal (for the Killing form) ideals:
\begin{equation}\label{dec-i}
\g=\oplus_{i=0}^{r}\g_i
\end{equation}
with $\g_0$ abelian and $(\g_i,\sigma_{\mid\g_i})$ irreducible, and let $\g_i=\hk_i\oplus\mk_i$ be the eigenspace decomposition of $\sigma_{\mid\g_i}$. Suppose now that we have an order four automorphism $\tau\colon\g\to\g$
such that $\tau^2=\sigma$.\\
There exists an unique decomposition of $\g$:
\begin{equation}\label{dec'-i}
\g=\g_0\oplus(\oplus_{i=1}^{r'}\g_i')
\end{equation}
where $\g_i'=\g_i$ or $\g_i'=\g_i\oplus\g_j$ with $(\g_i,\sigma_{\mid\g_i})$ and $(\g_j,\sigma_{\mid\g_j})$
isomorphic, such that $\tau(\g_i')=\g_i'$. Besides if $\g_i'=\hk_i'\oplus\mk_i'$ is the eigenspace decomposition of $\sigma_{\mid\g_i'}$, then $\tau(\hk_i')=\hk_i'$, $\tau(\mk_i')=\mk_i'$. Moreover if 
$\g_i'=\g_i\oplus\g_j$ then $\tau(\g_i)=\g_j$, $\tau(\hk_i)=\hk_j$, $\tau(\mk_i)=\mk_j$. Hence $\taum=\oplus_{i=0}^{r'}\tau_{\mk_i'}$ with $\tau_{\mk_i'}\in\Sigma(\tau_{\mk_i'})$, and $\tau_{\mid\g_i'}$ is the automorphism of $\g_i'$ defined by 
$ \tau_{\mk_i'}$ according to theorem~\ref{m} and we have $\tau_{\mid\g_i'}^2=\sigma_{\mid\g_i'}$.
\end{thm}
\begin{cory}\label{g-i'}
Let $M$ be a $G$-symmetric space, $G\subset \Is(M)$ and $\tau\colon\g\to\g$ an order four automorphism with $\tau^2=\sigma$. Let $M$ be its universal covering, which has a symmetric realisation $\tl M=\tl G/\tl H^0$, with $\pi\colon\tl G\to G$ a covering of $G$, $\tl H=\pi^{-1}(H)$ and $\tl G\subset\Is(\tl M)$, such that $\tau$ integrates into $\tl\tau\colon \tl G\to\tl G$. Then the decomposition of $\g$,(\ref{dec-i}), gives rise to the following decomposition of $\tl G$:
$$
\tl G=L_0\times L_1\times\cdots \times L_r
$$
 which induces the following decomposition of $\tl H^0$, corresponding also to the decomposition 
$\hk=\oplus_{i=0}^{r} \hk_i$:
$$
\tl H^0=H_0\times H_1\times\cdots\times H_r.
$$
Then $\tl\sigma$ admits the decomposition  $\tl\sigma=\prod_{i=0}^r \tl\sigma_i$ (with $\tl\sigma_i$ involution of $L_i$)
 and $H_i=(L_i^{\tl\sigma_i})^0$. Moreover there exists an unique decomposition of $\tl G$:
\begin{equation}\label{decgi'}
\tl G= L_0'\times L_1'\times\cdots\times L_{r'}'
\end{equation}
where $L_i'=L_i$ or $L_i'=L_i\times L_j$ with $(L_i,{\tl\sigma}_i)$ and $(L_j,{\tl\sigma_j})$ isomorphic.
Then $\tau$ admits the decomposition $\tl\tau=\prod_{i=0}^{r'}\tl\tau_i'$ with $\tl\tau_i'$ order four automorphism of $L_i'$.
Further, by identifying $(L_i,\tl\sigma_i)$ and $(L_j,\tl\sigma_j)$ (when $L_i'=L_i\times L_j$), then in (\ref{decgi'}),
we have either $L_i'=L_i$ and then $\tl\tau_i'=\tl\tau_i$ is an order four automorphism of $L_i$ so that 
$(L_i')^{\tl\tau_i'}=(L_i)^{\tl\tau_i}$, or $L_i'=L_i\times L_i$ and then 
$$
\tl\tau_i'\colon (a,b)\in L_i\times L_i\mapsto (\sigma_i(b), a)\in L_i\times L_i
$$
so that $(L_i')^{\tl\tau_i'}=\Delta(H_i)\subset H_i\times H_i$. Hence $\tl M=M_0\times M_1\times\cdots\times M_r$ 
with $M_i=L_i/H_i$, and $\tl N= N_0'\times N_1'\times\cdots\times N_{r'}'$ where either $N_i'=N_i=L_i/(L_i)_0$ is a 4-symmetric bundle over $M_i$, or $N_i'=L_i\times L_i/\Delta(H_i)$ is a 4-symmetric bundle over $M_i\times M_i=
L_i\times L_i/H_i \times H_i$ (and the fibre $H_i \times H_i/\Delta(H_i)\simeq H_i$ is a group).
\end{cory}
\textbf{Proofs of theorems~\ref{alpha},\ref{gi} and corollaries~\ref{alpha'},\ref{g-i'}} Use the fact that $\taum$ leaves invariant the metric in $\mk$ and the restriction to $\mk$ of the Killing form.\hfill$\blacksquare$
%%%%%%%%%%%%%%%%%%%%%%%%%%%%%%%%%%%%%%%%%%%%%%%%%%%%%%%%%%%%%%%%%%%%%%%%%%%%%%%%%%%%%%%%%%%%%%%%%%%%%%%%

\subsection{The semisimple case}\label{semisimple}
\begin{defn}
We will say that the Riemannian symmetric space $M$ is of semisimple type if $\Is(M)$ is semisimple.
\end{defn}
\begin{thm}
If $M$ is of semisimple type then each (connected) locally 4-symmetric bundle over $M$ is maximal and in the form $\und{N}=\Is(M)^0/\und{G_0}$, i.e. is an $\Is(M)^0$-orbit in 
$\Sigma(M)\cap \aut(M)$. In other words the set of locally 4-symmetric bundles over $M$ is $\mathcal N(M)=(\Sigma(M)\cap \aut(M))/\Is(M)^0\subset \Sigma(M)/\Is(M)^0$.
 \end{thm}
\begin{rmk}\emph{
The "size" of a maximal (locally) 4-symmetric bundle over $M$ in the twistor bundle $\Sigma(M)$ depends on the "size" of the isotropy subgroup $\iso(M)$ and on $J_0\in\Sigma(T_{p_0}M)$. In other words, if we want  a fibre $S_0\subset
\Sigma(T_{p_0}M)$ of maximal dimension, we must find $J_0\in\Sigma(T_{p_0}M)\cap\aut(T_{p_0}M)\supset
\Sigma(T_{p_0}M)\cap\iso(M)$  such that $T_{J_0}\und{S_0}=\und{\g_2}(J_0):=\{A\in\mathfrak{Is}_{p_0}(M)\mid AJ_0 + J_0A=0
\} $ is of maximal dimension, or equivalently such that $ \und{\g_0}(J_0)=\{ A\in\mathfrak{Is}_{p_0}(M)\mid AJ_0 - J_0A=0\}$ is of minimal dimension.
}\end{rmk}
\begin{rmk}\emph{
It is possible that there exist different non-isomorphic locally 4-symmetric bundles over $M$ (see section~\ref{cgrass}). And it is also possible that there does not exist any
locally 4-symmetric bundle over $M$. For example: $M=S^1\times S^3$, then $\Is(M)=SO(2)\times SO(4)$ and $\iso(M)=SO(3)$, and there does not exist $J_0\in\Sigma(\R^4)$ such that 
$J_0 SO(3)J_0^{-1}=SO(3)$.}
\end{rmk}
Moreover we have the following obvious theorem (see also \cite{jimenez}):
\begin{thm}
Let $(\g,\sigma)$ be an orthogonal symmetric Lie algebra. Then set $\g^*=\hk\oplus i\mk$ and $\sigma^*=\Id_{\hk}\oplus
-\Id_{i\mk}$. Then $(\g^*,\sigma^*)$ is an orthogonal symmetric Lie algebra. If $(\g,\sigma)$ is of the compact type then
$(\g^*,\sigma^*)$ is of the non-compact type and conversely. Now, for $\taum\in\mrm{End}(\mk)$, set
$\taum^*\colon iv\in i\mk\mapsto i\taum(v)$. Then 
$$
\taum\in\aut(\mk)\Longleftrightarrow \taum^*\in\aut(i\mk)
$$
 and $\taum\in\Sigma(\mk)$ \iif $\taum^*\in\Sigma(i\mk)$. In this case ($\taum\in\aut(\mk)\cap\Sigma(\mk)$) let $\tau$ (resp. $\tau^*$) be the automorphism of $\g$ (resp. $\g^*$) 
 defined by $\taum$ (resp. $\taum^*$) and denoting by $A^\C\in\mrm{End}(V^\C)$ the extension to $V^\C$ of $A\in\mrm{End}(V)$ ($V$ real vector space) then we have 
$$
\tau^\C={\tau^*}^\C\qquad\text{i.e. }\tau^*=\tau^\C_{\mid\g^*}
$$
\end{thm}
\begin{thm}
Let $M$ be an irreducible symmetric spaces of type II (compact type) or type IV (non-compact type) then there does not exist any (non-trivial) locally 4-symmetric bundle over $M$. Equivalently $\aut(M)\cap\Sigma(M)=\varnothing$, in other words, there does not exist any automorphism $\tau$ of $\mathfrak{Is}(M)$ such that $\tau^2=\sigma$.
\end{thm}
\textbf{Proof.} By duality, it is enough to prove the assertion for the compact type. In this case let $\tl M$ be the universal covering of $M$, we have $\tl M=H\times H/\Delta(H)$ and $\tl\sigma\colon (a,b)\in G\times G\mapsto (b,a)$.
Then an automorphism $\tau\colon\g\to \g$ must send $\g_1=\hk\oplus\{0\}$ either on $\g_1$ or on $\g_2=\{0\}\oplus\hk$ 
and idem for $\g_2$, and thus for any automorphism we have $\tau^2(\g_i)=\g_i$ and hence we cannot have $\tau^2=\sigma$.
This completes the proof.\hfill$\blacksquare$
%
%%%%%%%%%%%%%%%%%%%%%%%%%%%%%%%%%%%%%%%%%%%%%%%%%%%%%%%%%%%%%%%%%%%%%%%%%%%%%%%%%%%%%%%%
%
\subsection{The Euclidean case}\label{euclidean}
\begin{thm}
Let $M=\R^{2n}$ with its canonical inner product. Then $\Is(M)=O(2n)\ltimes\R^{2n}$ the group of affine isometries in $\R^{2n}$. Hence for any  $p_0\in\R^{2n}$, we have $\iso(M)=\{(F,(\Id-F)p_0), F\in O(2n)\}\simeq O(2n)$. In particular for $p_0=0$, $\iso(M)=O(2n)$. Thus we have $\forall p_0\in\R^{2n}$, $\Is(M)=\iso(M)\ltimes\R^{2n}$.\\
Further $M=G/H$ is a symmetric realisation with $G$ acting effectively \iif $G=H\ltimes\R^{2n}$ with $H\subset \iso(\R^{2n})$ for some $p_0\in \R^{2n}$. Then we have $G=H_0\ltimes\R^{2n}$ with $H_0=\mrm{pr}_{O(2n)}(H)\subset O(2n)$. The  involution for this realisation is 
$$ 
\sigma=\Int(-\Id,2p_0)\colon(h,x)\in G\mapsto (h,2(\Id-h)p_0- x)
$$
giving rise to the symmetry around $p_0$: $\sigma_0\colon x\in \R^{2n}\mapsto -(x-p_0)+p_0\in \R^{2n}$.
Let us fix $p_0=0$, so that for any symmetric realisation we have $H\subset \Is_{p_0}(M)= O(2n)$ and $\sigma=\Int(-\Id,0)$.\\
All (connected) locally 4-symmetric bundles over $M$ are globally 4-symmetric bundles over $M$. The twistor bundle, $\Sigma(\R^{2n})\times \R^{2n}$, is a globally 4-symmetric bundle over $M$. All the (connected) 4-symmetric bundles over $\R^{2n}$ are in the form: $S_0\times\R^{2n}$ where $S_0$ is a compact Riemannian symmetric space embedded\footnote{only immersed if $H$ is not closed in $O(2n)$} in $\Sigma^\eps(\R^{2n})$. Besides $\aut(T_{p_0}M)=\Is_{p_0}(M)=O(2n)$ so that any $J_0\in\Sigma(\R^{2n})$ defines the maximal 4-symmetric bundle $\Sigma(\R^{2n})\times \R^{2n}=(O(2n)\ltimes\R^{2n})/U(\R^{2n},J_0)$.
\end{thm}
\begin{thm}
Let $M$ be an Euclidean Riemannian symmetric space (i.e. its universal covering is an Euclidean space $\R^{2n}$). Then $M=\R^{2p}\oplus\mathbb{T}^{2q}$, $\Is(M)=O(2p)\times
(\mathfrak{S}_{2q}\otimes\{\pm 1\})\ltimes M$ ($\mathfrak{S}_{2q}$ is the group of permutations) and denoting by $\pi\colon\R^{2n}\to M$ the universal covering, and $p_0=\pi(0)$, 
then $\iso(M)=O(2p)\times(\mathfrak{S}_{2q}\otimes\{\pm 1\})$. 
Moreover $\aut(\R^{2n})=O(2n)$, and $J_0\in\Sigma(\R^{2n})$  defines the (connected) maximal 4-symmetric bundle over $M$: 
$(\Sigma(E^{2l})\times\{{J_0}_{\mid {E^{2l}}^\perp}\})\times M$, where $E^{2l}$ is the (unique) maximal subspace in $\R^{2p}$ invariant by $J_0$. In particular, $\aut^*(M)\cap
\Sigma(M)=\Sigma(\R^{2p})\times\Sigma(\R^{2q})\times M$.
\end{thm}
\textbf{Proof.} Let $\tl\pi\colon\tl G\to G$ be  a covering of $G=\Is(M)^0$ acting symmetrically and effectively on $\tl M=\R^{2n}$ and $\tl\sigma\colon\tl G\to\tl G$ the 
corresponding involution. Then setting $\tl H=(\tl G^{\tl\sigma})^0$, we have according to the previous theorem $\tl G=\tl H\ltimes\R^{2n}$ and $\tl H\subset SO(2n)$. Then setting 
$D=\ker\pi$, $D$ is a discrete central subgroup of $\tl G$. Besides it is easy to see that $\mrm{Cent}(\tl G)= \mrm{Cent}(\tl H\times \R^{2n})
=\R^{2q}$ where $\R^{2q}$ is the maximal subspace of $\R^{2n}$ fixed by $\tl H$, i.e. $\tl H\subset SO(2p)\times\{\Id_{2q}\}$ ($2p+2q=2n$). Hence $D=\oplus_{i=1}^r\Z e_i$ with 
$(e_i)_{1\leq i\leq r}$ $\R$-free so that $G=\Is(M)^0=\tl G/D=\tl H\ltimes M'$ with $M'=\R^{2p}\oplus\R^{2q-r}\oplus\R^{r}/\mathbb{Z}^{r}$. Moreover we have $\sigma\colon (h,x)
\in\tl H\times M'\to (h,-x)$ because $\tl \sigma=\Int(-\Id,0)$ (see the previous theorem) and thus $G^\sigma=\tl H$ but the isotropy subgroup of $G$ at $p_0$ satisfies 
$H\supset\tl\pi(\tl H)$ (because $\tl H$ is connected), but $\tl\pi(\tl H)=\tl H$ ($D\cap\tl H=\{1\}$) and thus $H=\tl H$. Thus $M=G/H=M'$. Now, we have to compute $\Is(M)$, 
we know that $\Is(M)^0=H\ltimes M\subset SO(2p)\ltimes M$. In the other hand, any $g\in\Is(M)$ can be lifted into $\tl g\in O(2n)\ltimes\R^{2n}$, and conversely $\tl g\in O(2n)\ltimes\R^{2n}$ corresponds to some $g\in\Is(M)$ \iif $\tl g(D)=D$ which is equivalent to $\tl g\in[O(2p+2q-r)\times (GL_r(\Z)\cap O(\R^r))]\ltimes \R^{2n}=[O(2p+2q-r)\times (\mathfrak{S}_r\ltimes\{\pm\Id\})]\ltimes \R^{2n}$. Hence $\iso(M)^0=SO(2p+2q-r)$ and thus $r=2q$. Finally $M=\R^{2p}\oplus\mathbb{T}^{2q}$, $\Is(M)=O(2p)\times (\mathfrak{S}_r\ltimes\{\pm\Id\})\ltimes M$,  
$\iso(M)=O(2p)\times (\mathfrak{S}_r\ltimes\{\pm\Id\})$, and $\iso(M)^0=H=SO(2p)$. We conclude by remarking that $J_0\in \Sigma(\R^{2n})$ satisfies $J_0HJ_0^{-1}=H$  for $H\subset
SO(2p)$ connected and maximal \iif $H=SO(E^{2l})$ and $J_0\in\Sigma(E^{2l})\times\Sigma({E^{2l}}^\perp)$. This completes the proof.\hfill$\blacksquare$
\begin{rmk}\emph{
We can use the second elliptic integrable system in the Euclidean case to "modelize" this system in the general case. Indeed, let us consider $M$ a Riemannian symmetric space of the semisimple type (then its isotropy subgroup $H=\iso(M)$ is essentially its holonomy group, i.e. they have the same identity component) with $\tau\colon\g\to\g$ an order four automorphism such that $\tau^2=\sigma$. Then we can associate to the corresponding locally 4-symmetric bundle $N$ over $M$, the 4-symmetric bundle over $M_0=\mk=H\ltimes \mk/H$ : $N_0=H\ltimes \mk/G_0=S_0\times\mk\subset\Sigma(\mk)\times \mk$, and to the second elliptic integrable system in $N$, its "linearized" in $N_0$. We conjecture that the "concrete" geometrical interpretation (i.e. in terms of the second fundamental form of the surface $X$ etc...)
 is the same for the linearized and the initial system. This is what happens in dimension 4.
}\end{rmk}
\begin{rmk}\emph{
The second elliptic integrable system can be viewed as "a couplage" between the harmonic map equation in $S_0=H/G_0$ and a kind of  Dirac equation in $\g_{-1}$: $\partial_{\bar z}u_1 + [\bar u_0,u_1] + [\bar u_1, u_2] =0$. In the Euclidean case, the projection on the "group part", $\g=\hk\ltimes \mk\to \hk$, of the second elliptic system  is only the harmonic map equation in $H/G_0$. In other words, the second elliptic integrable system is only the harmonic map equation in $H/G_0$ and a kind of Dirac equation in $\C^n$ ($\cong(\g_{-1},J_0)$). In particular, if we apply any method of integrable systems theory using loop groups (DPW, Dressing action etc..) or something else (spectral curves) to the second elliptic  system in $G/G_0$ and then project in the group part ($\mrm{pr}\colon H\ltimes\mk\to H$), we obtain the same method applied to the first elliptic integrable system in $H/G_0$ i.e. the harmonic map equation in $H/G_0$. For example, if we apply the DPW method: given $\mu=(\mu_\hk,\mu_\mk)$ a holomorphic potential, we have $\mrm{pr}(\mathcal{W}_{G/G_0}(\mu))= \mathcal{W}_{H/G_0}(\mu_\hk)$ where $\mathcal{W}_{G/G_0}, \mathcal{W}_{H/G_0}$ are the Weierstrass representations for each elliptic system. Hence to solve the second elliptic system, we can first solve the harmonic map equation in $H/G_0$, by using  any method of integrable systems theory which gives us a lift $h$ in $H$ of a harmonic map in $H/G_0$, and then we have to solve the Dirac equation with parameters $u_0,u_2$ given by the lift : $h^{-1}\partial_{z}h=u_0 + u_2 $ following $\hk=\g_0\oplus\g_2$ (see \cite{ki1}). However, the Dirac equation is not intrinsic since it depends on the lift $h$ of the harmonic map (see \cite{ki1}).\\
In the particular case where $S_0$ is a group and $H=G_0\rtimes S_0$, (for example $S_0=G_0\times G_0/G_0$), then we have a canonical lift and then the Dirac equation becomes intrinsic (see \cite{ki1}). It is in particular what happens for Hamiltonian stationary Lagrangian surfaces : in $\C^2$ we have an intrinsic Dirac equation whereas in the others Hermitian symmetric spaces this equation does not exist (see \cite{HR1,HR2,HR3}). It is also what happens in \cite{ki1} when we take  for $S_0$ the subsphere $S^3\subset S^6$ ($S^6$\ embeds in $\Sigma^+(\R^8)$ by the left multiplication in $\oct$).
}\end{rmk}
%
%%%%%%%%%%%%%%%%%%%%%%%%%%%%%%%%%%%%%%%%%%%%%%%%%%%%%%%%%%%%%%%%%%%%%%%%%%%%%%%%%%%%%%%%
%
\section{Examples of 4-symmetric bundles}\label{example}
We use the notations of section~\ref{4-sym}.
\subsection{The sphere} 
Let us consider $M=S^{2n}=SO(2n+1)/SO(2n)$ with $G=SO(2n+1)$, $H=SO(2n)$ and the involution $\sigma=\Int(\mrm{diag}(\Id_{2n},-1))$. Then $G^{\sigma}=SO(2n)\bigsqcup O^-(2n)\times\{-1\}$. Hence $H=(H^\sigma)^0$, $M_{min}=\R\mathbb{P}^{2n}$ and 
$M_{max}=S^{2n}$.\footnote{$M_{max}$ is simply connected and  $M_{min}$ is  the adjoint space.} We have also 
$$
\hk=\mathfrak{so}(2n),\qquad \mk=\left\{\begin{pmatrix} 0 & v\\ -v^t & 0 \end{pmatrix}, v\in \R^{2n}\right\}=
\{i_{\mk}(v), v\in\R^{2n}\}
$$
where $i_{\mk}\colon\R^{2n}\to\mk$ is defined in an obvious way. Now, let us consider the action of $H$ on $\mk$:
for $h\in SO(2n)$, $\xi=i_{\mk}(v)\in\mk$, we have
$$
\Adm h(\xi)=i_{\mk}(h.v)
$$
hence  $K=\ker\Adm=\{\Id\}$ and the action of $G$ is effective (in fact $SO(2n+1)$ is simple because $2n+1$ is odd).
Identifying $\mk$ with $\R^{2n}$ via $i_{\mk}$ we have: $\forall h\in SO(2n), \Adm h=h$ i.e. $\Adm=\Id$.
Moreover $SO(2n+1)$ is the connected isometry group of $S^{2n}$. Now, according to theorem~\ref{m}, define a locally 
4-symmetric bundle over $M=S^{2n}$ is equivalent to give ourself $\taum\in\Sigma(\mk)\cap\aut(\mk)=\Sigma(\mk)$. Further, given $J_0\in\Sigma^\eps(\R^{2n})$, let us define the order four automorphism of $G$: $\tau=\Int(\diag(-J_0,1))$. Then 
$\tau^2=\sigma$ and since $\tau_H=\Int J_0$ and $\tm=J_0$, we obtain all the locally 4-symmetric bundles over $M$ which are all globally 4-symmetric bundles over $M$.\\
Moreover, we have $G^{\tau}=\mrm{Com}(J_0)\cap SO(2n)=U(\R^{2n},J_0)$. Hence $G^{\tau}=(G^{\tau})^0=G_0$ thus $S_0=H/G_0=
\Int(SO(2n))(J_0)=\Sigma^\eps(\R^{2n})$ and thus $N=G/G_0=\Sigma^\eps(S^{2n})$.
\subsection{Real Grassmannian}
More generally, let $p,q\in\mathbb{N}^*$ such that $pq$ is even and let us consider $M=SO(p+q)/SO(p)\times SO(q)=Gr_p(\R^{p+q})$ (oriented $p$-planes in $\R^{p+q}$). Since $p$ and $q$ play symmetric roles, we will suppose that $p$ is even and that it has the biggest divisor in the form $2^r$.  We have $\dim M=pq$ and the following setting
$$
\begin{array}{l}
G=SO(p+q),\ H=SO(p)\times SO(q);\ \sigma=\Int(\diag(\Id_p,-\Id_q))\text{ and}\\
G^{\sigma}=SO(p)\times SO(q)  \bigsqcup O^-(p)\times O^-(q).
\end{array}
$$
Then $H=(G^\sigma)^0$ so that $M_{min}=Gr_p^*(\R^{p+q})$ (non-oriented $p$-planes in $\R^{p+q}$) and $M_{max}=Gr_p(\R^{p+q})=M$. Besides $\hk=\so(p)\oplus\so(q)$, and $\mk=\left\{\begin{pmatrix} 0 & B\\ - B^t & 0 \end{pmatrix}, B\in\mak{gl}_{p,q}(\R)
\right\}=i_{\mk}(\mak{gl}_{p,q}(\R))$ ($i_{\mk}$ defined in an obvious way).\\
Now let us compute $\Adm$. For $h=\diag(A,C)$ and $\xi=i_{\mk}(B)$, we have:
$$
\Adm h(\xi)=i_{\mk}(ABC^{-1}).
$$
Under the identification $i_{\mk}$ we have $\Adm(A,C)=L_AR_{C^{-1}}=\chi(A,C)$, by introducing the morphism $\chi\colon (A,C)\in GL_p(\R)\times GL_q(\R)\mapsto L(A)R(C^{-1})\in GL(\mak{gl}_{p,q}(\R))$. Hence $K=\ker\Adm=\{\pm\Id\}$ if $q$ is even and $K=\{\Id\}$ if not. Thus the connected isometry group of $M$, $\Is(M)^0$, is $G'=G/K=PSO(p+q)$ if $q$ is even and $G'=G=SO(p+q)$ if not. Let us compute $\aut(\mk)$: we already know that $\aut(\mk)\supset H\supset \aut(\mk)^0$. But, it is well known that the automorphisms of $\so(n+1)$ are all inner automorphisms by $O(n+1)$ so we have $\aut(\mk)=\{L_AR_{C^{-1}},(A,C)\in O(p)\times O(q)\}$. Thus $J_0=L(J_1)R(J_2^{-1})\in \aut(\mk)$ is in $\Sigma(\mk)$ \iif:
$$
\left\{\begin{array}{lccl}
(J_1^2,J_2^2) & = & \pm (-\Id_p,\Id_q) & \text{ if } $q$ \text{ is even,}\\
(J_1^2,J_2^2) & = &  (-\Id_p,\Id_q) & \text{ if } $q$ \text{ is odd.}
\end{array}\right.
$$
Then the associated order four automorphism is $\tau=\Int(\diag(J_1,J_2))$. In particular, $\tau(H)=H$ and $\tau_H=\Int 
J_1\times \Int J_2$. Besides, $\aut(\mk)\cap\Sigma(\mk)$ has respectively $2(p+q+2)$ or $2(q+1)$ connected components if $q$ is even or $q$ is odd respectively. Each connected component is an $\Adm H$-orbit and corresponds to the fibre of a different maximal 4-symmetric bundle over $M$.\\
 Moreover to fix ideas let us suppose that we have $J_1\in \Sigma(\R^p),J_2\in OS(\R^q)$, the set of orthogonal symmetries in $\R^q$, then $G^\tau=U(\R^p,J_1)\times S(O(E_1)\times O(E_2))$ with $E_1=\ker(J_2-\Id)$, $E_2=\ker(J_2+\Id)$. We have 
$G^\tau\subset H$. Let $OS_r(\R^q)=\Int(SO(q))(\Id_r,-\Id_{q-r})$ be the set of orthogonal symmetries in $\R^q$ with $\dim E_1=r$. Then $H/G^\tau=\Int(H)(J_1,J_2)=\Sigma^\eps(\R^p)\times OS_r(\R^q)$ ($\eps$ being determined by $J_1$) and 
\begin{equation}\label{not/4}
G/G^\tau=\{(x,J),x\in M,J\in\Sigma^\eps(x)\times OS_r(x^\perp)\}.
\end{equation}
Now let us compute $G_0$ according to (\ref{g-0}): $h=(A,C)\in H$ is in $G_0$ \iif $\Adm\tau(h)=\Adm h$ i.e.: if $q$ is odd, $\tau(h)=h$, and $G_0=G^\tau\cap H=G^\tau$; if $q$ is even, $\tau(h)=\pm h$ (and $G_0=\pi_K^{-1}(G_0')$ with $G_0'={G'}^{\tau'}\cap H'$), i.e. $h\in G^\tau$ or $\tau(h)=-h$. The existence of solutions of this last equation depends on $p,q$ and $r$ (we remark that if $h_1$ is a solution then the set of solutions is $h_1G^\tau$). One finds that the equation $\tau(h)=-h$ ($q$ is even) has a solution in $G^\sigma$ \iif $\dim E_1=\dim E_2=q/2$ and that this solution is in $H$ if $p/2$ is even and in $O^-(p)\times O^-(q)$ (the other component of $G^\sigma$)  if $p/2$ is odd. Hence, if $p$ is divisible by 4, $q$ is even and $r=q/2$ (i.e. $J_0\in\chi(\Sigma(\R^p)\times OS_{q/2}(\R^q))$), we have $G_0=G^\tau\bigsqcup h_1 G^\tau$. In all the other cases we have $G_0=G^\tau$.\\
In conclusion, let us denote by $N^L(r,\eps):=N(J_0)$ (resp. $N^R(r,\eps)$)  the maximal 4-symmetric bundle over $M$ corresponding to $J_0\in\chi(\Sigma^\eps(\R^p)\times OS_r(\R^q))$ (resp. $\chi(OS_r(\R^p)\times\Sigma^\eps(\R^q)$). Then:
\begin{description}
\item if $p$ is not divisible by 4 or $q$ is odd, $N^\alpha(r,\eps)$ is given by (\ref{not/4}), for all $(\alpha,r,\eps)$, 
\item if $p$ is divisible by 4, $q$ even not divisible by 4 then for $(\alpha,r)\neq (L,q/2)$, $N^\alpha(r,\eps)$ is given by (\ref{not/4}) and for $(\alpha,r)= (L,q/2)$ it is given by (\ref{/4}), below,
\item if $p$ and $q$ are divisible by 4, then for $(\alpha,r)\in \{(L,q/2),(R,p/2)\}$, $N^\alpha(r,\eps)$ is given by (\ref{/4}), and for the other choices it is given by (\ref{not/4}),
\end{description} 
\begin{equation}\label{/4}
\begin{array}{rcl}
N^L(r,\eps) & = & \{(x,J),x\in M,J\in P(\Sigma^\eps(x)\times OS_r(x^\perp))\} \\
N^R(r,\eps) & = & \{(x,J),x\in M,J\in P(OS_r(x)\times \Sigma^\eps(x^\perp))\}
\end{array}
\end {equation}
where $P(\Sigma^\eps(x)\times OS_r(x^\perp))=\Sigma^\eps(x)\times OS_r(x^\perp)/\{\pm\Id\}$. In the cases described by (\ref{/4}), $G/G^\tau$ is not a submanifold of $\Sigma(M)$.
%
%%%%%%%%%%%%%%%%%%%%%%%%%%%%%%%%%
%
\subsection{Complex Grassmannian}\label{cgrass}
Let us consider $M=SU(p+q)/S(U(p)\times U(q))=Gr_{p,\C}(\C^{p+q})$. We have $\dim M=2pq$ and the following setting
$$
\begin{array}{l}
G=SU(p+q),\ H=S(U(p)\times U(q));\ \sigma=\Int(\diag(\Id_p,-\Id_q))\text{ and}\\
G^{\sigma}=H=(G^{\sigma})^0.
\end{array}
$$
Besides $\hk=\mak{s}(\mak{u}(p)\oplus\mak{u}(q))$ and $\mk=\left\{\begin{pmatrix} 0 & B\\ - B^* & 0 \end{pmatrix}, B\in\mak{gl}_{p,q}(\C)\right\}=i_{\mk}(\mak{gl}_{p,q}(\C))$. Let us compute $\Adm$. For $h=\diag(A,C)$ and $\xi=i_{\mk}(B)$, we have:
$$
\Adm h(\xi)=i_{\mk}(ABC^{-1}).
$$
Under the identification $i_{\mk}$ we have $\Adm(A,C)=L_AR_{C^{-1}}=\chi(A,C)$, by introducing the morphism $\chi\colon (A,C)\in GL_p(\C)\times GL_q(\C)\mapsto L(A)R(C^{-1})\in GL(\mak{gl}_{p,q}(\C))$\footnote{For the following it useful to keep in mind that we have $\Adm H=\chi(S(U(p)\times U(q)))=\chi(U(p)\times U(q))$ and $\ker\chi=\C^*\Id$.
}. Hence $K=\ker\Adm=\{(\lm\Id_p,\lm\Id_q), \lm\in\C,\lm^{p+q}=1\}=\Hat{U}_{p+q}\Id\simeq\Z_{p+q}$ (with $\Hat{U}_{p+q}=\exp(\frac{2i\pi}{p+q}\Z)$). Thus $G'=G/K=PSU(p+q)$ and $H'=S(U(p)\times U(q))/\Hat{U}_{p+q}\simeq S(U(p)\times U(q))$. The connected isometry group is the unitary group of $M$: $\Is(M)^0=U(M)=G'=PSU(p+q)$.\\
It is well known that the group of automorphisms of $SU(p+q)$ has two components (the $\C$-linear one and the anti-$\C$-linear one) and is generated by the inner automorphisms and the complex conjugation: $g\in SU(p+q)\mapsto \bar g\in SU(p+q)$. In particular, $\aut(\mk)=\Adm H\rtimes\{\Id,c\}=\chi(S(U(p)\times U(q))\cdot\{(\Id,\Id),(b_p,b_q)\})$ with $c=L(b_p)R(b_q^{-1}):B\in\mak{gl}_{p,q}(\C)\mapsto \Bar B\in \mak{gl}_{p,q}(\C)$, $b_n\colon v\in\C^n\mapsto\bar v\in \C^n$.\\
The complex structure in $\mk=\mak{gl}_{p,q}(\C)$ is defined by $L(I_p)=R(I_q)$ where $I_n=i\Id_n$ is the canonical complex structure in $\C^n$, and the two connected components of $\aut(\mk)$ are respectively the elements in $\aut(\mk)$ which commute and those which anticommute with this complex structure.\\
Moreover, $J_0=L(J_1)R(J_2^{-1})\in\aut(\mk)^0=\Adm H$ is in $\Sigma(\mk)$ \iif $(J_1^2,J_2^2)\in(-\Id_p,\Id_q)U(1)$. Then let us set $\Sigma_\lm=\{(J_1,J_2)\in U(p)\times U(q)\mid (J_1^2,J_2^2)=\lm (-\Id_p,\Id_q)\}$. Then we have $\chi(\Sigma_\lm)=\chi(\Sigma_0)$ for all $\lm\in U(1)$ since $\Sigma_\lm=\lm^{\frac{1}{2}}\Sigma_0$ with $\lm^{\frac{1}{2}}$ a root of $\lm$. Thus according to the following lemma, $\aut(\mk)^0\cap\Sigma(\mk)$ has $(p+1)(q+1)$ connected components, which are
$\Adm H$-orbits and correspond to the fibres of  different maximal 4-symmetric bundles over $M$.
\begin{lemma} 
Let $J\in U(n)$, then $J^2=-\Id$ (resp. $J^2=\Id$) \iif there exists $h\in U(n)$ such that $hJh^{-1}=\diag(i\Id_l,-i\Id_{n-l})$ for some $l\in\{0,\ldots,n\}$ (resp. $hJh^{-1}=\diag(\Id_r,-\Id_{n-r})$ for some $r\in\{0,\ldots,n\}$).
\end{lemma}
Then the order four automorphism corresponding to $J_0$ is $\tau=\Int(\diag(J_1,J_2))$, with\footnote{$\mrm{I}_{l,p-l}=\diag(\Id_l,-\Id_{p-l})$} $J_1\in\Ad U(p)(i\mrm{I}_{l,p-l})\cong iGr_{l,\C}(\C^p)$, $J_2\in \Ad U(q)(\mrm{I}_{r,q-r})\cong Gr_{r,\C}(\C^q)$. Hence $G^\tau=S(U(l)\times U(p-l)\times U(r)\times U(q-r))$; the fibre of the 4-symmetric space $G/G^\tau$ is $H/G^\tau=Gr_{l,\C}(\C^p)\times Gr_{r,\C}(\C^q)$, and 
\begin{equation}\label{g/gtau}
G/G^\tau=\{(x,J), x\in Gr_{p,\C}(\C^{p+q}), J\in Gr_{l,\C}(x)\times Gr_{r,\C}(x^\perp)\}.
\end{equation}
Further, $G_0$ is defined by: $\Adm\tau(h)=h,h\in H$, i.e. $(J_1AJ_1^{-1},J_2CJ_2^{-1})=\lm(A,C)$ for some $\lm\in K$. But it is easy to see that we must have $\lm^2=1$ and thus $\tau(h)=\pm h$. One finds that $\tau(h)=- h$ has solutions \iif $p,q$ are even and $l=p/2,r=q/2$. Finally, in the $\C$-linear case, the maximal 4-symmetric bundle $N=G/G_0$ is given by  
\begin{equation}\label{g/g0}
G/G_0=\{(x,J), x\in Gr_{p,\C}(\C^{p+q}), J\in Gr_{l,\C}(x)\times Gr_{r,\C}(x^\perp)/\Z_2\}
\end{equation}
if $p,q$ are even and $l=p/2,r=q/2$, and by (\ref{g/gtau}) in all the other cases.\\[0,1cm]
In the antilinear case,  $J_0=L(J_1)R(J_2^{-1})\in\aut(\mk)^0.c$, with $(J_1,J_2)=(J_1'b_p,J_2'b_q)$, is in $\Sigma(\mk)$ \iif $(J_1^2,J_2^2)=(J_1'\overline{J_1'},J_2'\overline{J_2'})\in (-\Id_p,\Id_q).U(1)$. It is easy to see that we can only have
\begin{equation}\label{j1j2}
(J_1^2,J_2^2)=\pm(-\Id_p,\Id_q).
\end{equation}
Hence according to the following lemma:\\
-- if $p,q$ are odd then $\Sigma(\mk)\cap(\aut(\mk)^0.c)=\varnothing$,\\
-- if $p,q$ are even then the two signs $\pm$ are realized in (\ref{j1j2}) and thus $\Sigma(\mk)\cap(\aut(\mk)^0.c)$ has 2 connected components,\\
-- if $p,q$ have opposite parities, then only one sign is realized in (\ref{j1j2}) and $\Sigma(\mk)\cap(\aut(\mk)^0.c)$ has one component.
\begin{lemma}
Let $E\subset \C^n$ be a Lagrangian n-plan, i.e. $E\overset{\perp}{\oplus} iE=\C^n$ and let $b_E$ be the associated conjugation: $v+iw\mapsto v-iw$ for $v,w\in E$. Then $U(n).b_E=b_E.U(n)$ does not depend on $E$ and is the set of anti-$\C$-linear isometries in $\C^n$ (the elements in $O(\R^{2n})$ which anticommute with the complex structure $I=i\Id$). Moreover for any $J$ in this set there exists a Lagrangian n-plane $E$ such that $J=J_E.b_E=b_E.J_E$ with $J_E\in O(E)$. Besides $J\in \Sigma(\R^{2n})$ (resp. $OS(\R^{2n})$) \iif $J_E\in\Sigma(E)$ (resp. $OS(E)$). In particular $\Sigma(\R^{2n})\cap (U(n).b_E)\neq \varnothing$ only if $n$ is even, moreover $\Sigma(\R^{2n})\cap (U(n).b_E)\subset \Sigma^+(\R^{2n})$. Then given any $J_1\in\Sigma(\R^n)$ (resp. $OS(\R^n)$) there exists $h\in U(n)$ such that $h.E=\R^n$, $hJ_Eh^{-1}=J_1$ and thus $hJh^{-1}=J_1.b_{\R^n}$.
\end{lemma}
Then the order four automorphism corresponding to $J_0$ is $\tau=\Int(\diag(J_1,J_2))$ with $J_1\in\Ad U(p)(J_{\frac{p}{2}}.b_p)$, $J_2\in \Ad U(q)(b_q)$ and $J_{\frac{p}{2}}=\begin{pmatrix} 0 & \Id_{\frac{p}{2}}\\ -\Id_{\frac{p}{2}} & 0
\end{pmatrix}$. In other words $J_1$ is any complex structure in $\R^{2p}$ anticommuting with $I_p$ and $J_2$ is any orthogonal conjugation in $\C^q$. Hence, we have $G^\tau=Sp(p/2)\times SO(q)$. Hence $U(p)\times U(q)/G^\tau=\Sigma^+(\C^p)_-\times\lag(\C^q)$ where $\Sigma^+(\C^p)_-=\Sigma(\R^{2p})\cap\mrm{Ant}(I_p)$ are the complex structures in $\R^{2p}$ anticommuting with $I_p$ and $\lag(\C^q)$ are the oriented Lagrangian planes in $\C^q$. Thus we have:
$$
H/G^\tau=S(\Sigma^+(\C^p)_-\times\lag(\C^q)):=\{(J,P)\in \Sigma^+(\C^p)_-\times\lag(\C^q)\mid {\det}_{\C}(J){\det}_\C(P)=1\}.
$$ 
It is easy to define ${\det}_{\C}$ on $\lag(\C^q)$; and for $\Sigma^+(\C^p)_-$, we set ${\det}_{\C}(J)={\det}_{\C}(E)$ for $E$ any Lagrangian $n$-plane invariant by $J$ (definition independent on the choice of $E$). Then
$$
G/G^\tau=\{(x,J,P), x\in Gr_{p,\C}(\C^{p+q}), (J,P)\in \Sigma^+(x)_-\times\lag(x^\perp)\}.
$$
Let us compute $G/G_0$. We have to solve for $(A,C)\in U(p)\times U(q)$: $(J_{\frac{p}{2}}\Bar A J_{\frac{p}{2}}^{-1}, \Bar C)=\lm (A,C)$ for $\lm\in U(1)$ whose the solutions are $\pm \lm^{\frac{1}{2}}(Sp(p/2)\times O(q))$. Hence we have $G_0'=G_0/K=\chi(U(1)(Sp(p/2)\times O(q)))=\chi(Sp(p/2)\times O(q))=$
$$
\begin{cases}
\chi(Sp(p/2)\times SO(q)) & \text{ if } q\text{ is odd}\\
\chi(G^\tau)\bigsqcup h_1\chi(G^\tau) & \text{ if } q\text{ is even.}
\end{cases}
$$
Then $G'/G_0'=G/G_0=U(p+q)/(U(1)(Sp(p/2)\times O(q)))=PSU(p+q)/P(Sp(p/2)\times O(q))$ hence $N=G/G_0$ is equal to $(G/G^\tau)/\Z_{p+q}$ if $q$ is odd and to $(G/G^\tau)/\Z_{2(p+q)}$ if $q$ is even.
%
%
%%%%%%%%%%%%%%%%%%%%%%%%%%%%%%%%%%%%%%%%%%%%%%%%%%%%%%%%%%%%%%%%%%%%%%%%%%%%%%%%%%%%%%%%
%
\section{Appendix}
\begin{thm}\label{compact}
Let $G$ be a connected Lie group with an involution $\sigma$. If $\Adm(G^\sigma)^0$ is compact  (resp. relatively compact) then $\Adm H$ is compact (resp. relatively compact) for any $H$ such that $(G^\sigma)^0\subset H\subset G^\sigma$.
\end{thm}
\textbf{Proof.} According to \cite{an} (lemma 2.7), $(G^\sigma)/(G^\sigma)^0$ is finite hence $H/(G^\sigma)^0$ is finite and the theorem follows.\hfill$\blacksquare$
\begin{cory}\label{cpair}
We give ourself the same setting and notations as in remark \ref{univcov}.\\
If $\tl H=(\tl G^{\tl\sigma})^0$ satisfies: $\Adm \tl H$ is compact (resp. relatively compact), then for any symmetric pair $(G,H)$, $\Adm H$ is compact (resp. relatively compact). In other words if one symmetric pair (associated to $(\g,\sigma)$) is Riemannian then all the others are also.
\end{cory}
\textbf{Proof.} Since $\tl G$ is simply connected, it is the universal covering of $G$ and we have a covering $\pi\colon\tl G\to G$. Then $\Adm\tl H=\Adm H^0$ (there are connected with the same Lie algebra) hence $\Adm H^0$ is compact and then according to the previous theorem, $\Adm H$ is compact.\hfill$\blacksquare$
\begin{cory}\label{autg}
Let $(G,H)$ be a symmetric pair with involution $\sigma$ and $\tau\colon G\to G$ an order four automorphism such that $\tau^2=\sigma$. Then if $\Adm H$ is compact (resp. relatively compact) then the subgroup generated by $\Adm H$ and $\tm$,
$\mrm{Gr}(\Adm H,\tm)$ is compact (resp. relatively compact).
\end{cory}
\textbf{Proof.}
We have $\tm(\Adm G^\sigma)\tm^{-1}= \Adm \tau(G^\sigma)= \Adm G^\sigma$. Hence $\mrm{Gr}(\Adm G^\sigma,\tm)=
(\Adm G^\sigma)\mrm{Gr}(\tm)$ which is (relatively) compact because so is $\Adm G^\sigma$, according to theorem~\ref{compact}, and then $\mrm{Gr}(\Adm G^\sigma,\tm)$ is (relatively) compact because since $\mrm{Gr}(\Adm H,\tm)
\supset (\Adm H)\mrm{Gr}(\tm)$  then $\Adm G^\sigma/\Adm H$ is a covering of $\mrm{Gr}(\Adm G^\sigma,\tm)/\mrm{Gr}(\Adm H,\tm)$ which is consequently finite.\hfill$\blacksquare$
\begin{thm}\label{g'}
Let $(G,H)$ be a symmetric pair with involution $\sigma \colon G\to G$ and $\tau\colon \g\to \g$ an order four automorphism such that $\tau^2=\sigma$. Then if $\Adm H$ is  relatively compact then the subgroup generated by $\Adm H$ and $\tm$,
$\mrm{Gr}(\Adm H,\tm)$ is  relatively compact.
\end{thm}
\textbf{Proof.} Let $G'=\Ad G$, then $C:=\ker\Ad=$\,center of $G$ and we can identify $\Ad$ to the covering $\pi\colon G\to G/C$ and $G'$ to $G/C$. The automorphism $\sigma$ gives rise to $\sigma'\colon G'\to G'$ such that $\sigma'\circ\pi=\pi\circ\sigma$. Besides the automorphism $\tau$ integrates in $G'$ into $\tau'$ defined by $\tau'=\Int\tau\colon\Ad g\in G'\mapsto\tau\circ\Ad g\circ\tau^{-1}$ and we have $\tau'\circ\pi=\pi\circ\tau$ and  $\tau'^2=\sigma'$. Then according to corollary~\ref{autg}, $\mrm{Gr}(\Adm {G'}^{\sigma'},\tm)$ is relatively compact since according to corollary~\ref{cpair}, 
$\Adm {G'}^{\sigma'}$ is relatively compact because $\Adm H$ is so. Moreover we have ${G'}^{\sigma'}\supset \pi(G^\sigma)$ then (since $\Ad\pi(g)=\Ad g\ \forall g\in G$) $\Adm {G'}^{\sigma'}\supset \Adm G^\sigma\supset\Adm H$ thus $\mrm{Gr}(\Adm H,\tm)$ is relatively compact.\hfill$\blacksquare$
\begin{thm}
Let $(\g,\sigma)$ be an orthogonal symmetric Lie algebra\footnote{i.e. $\sigma$ is an involutive automorphism and $\hk=\g^{\sigma}$ is compactly embedded in $\g$ (see \cite{Hel})} such that $\hk=\g^\sigma$ contains no ideal $\neq 0$ in $\g$.
Then for any  symmetric pair $(G,H)$ associated with $(\g,\hk)$, the associated symmetric space $M=G/H$ is Riemannian.
Moreover let $\tl G$ be the simply connected Lie group with Lie algebra $\g$, $\tl\sigma$ integrating $\sigma$, $\tl H=
(\tl G^{\tl\sigma})^0$ and $\tl C$ the center of $\tl G$. Then we have $\tl H=\tl G^{\tl\sigma}$. Further, for any subgroup $S$ of $\tl C$ put 
$$
H_{S}=\{g\in\tl G\mid \tl\sigma(g)\in g.S\}.
$$
The symmetric spaces $M$ associated with $(\g,\sigma)$ (i.e. $(G,H)$ is associated with $(\g,\hk)$) are exactly the spaces $M=G/H$ with
\begin{equation}\label{S}
G=\tl G/S\quad\mrm{ and }\quad H=H^*/S
\end{equation} where $S$ varies through all $\tl{\sigma}$-invariant subgroups of $\tl C$ and $H^*$ varies through all 
$\tl{\sigma}$-invariant subgroups of  $\tl G$ such that $\tl H S\subset H^*\subset H_S$. Hence, all the symmetric spaces 
$M=G/H=\tl G/H^*$ associated with $(\g,\sigma)$ cover the adjoint space of $(\g,\sigma)$: $M'=G'/{G'}^{\sigma'}=\tl G/ H_{\tl C}$\footnote{ with the notation of the proof of theorem~\ref{g'}. For any $(G,H)$ symmetric pair associated with $(\g,\sigma)$, we have $G'=\Ad G=\Int(\g)$ the group of inner automorphism of $\g$ (see \cite{Hel}) and $\sigma$ induces an automorphism $\sigma'$ of $G'=\Int(\g)$.} and are covered by $\tl M=\tl G/\tl H$ (the universal covering):
\begin{equation}\label{covering}
\tl M \rightarrow M \rightarrow M'.
\end{equation}
Besides if $\langle \cdot,\cdot\rangle$ is an $\Adm {G'}^{\sigma'}$-invariant inner product then it is invariant by $\adm H=\Adm H^*$ for any $H$ described above, and the coverings (\ref{covering}) are Riemannian, when $M,\tl M ,M'$ are endowed with the corresponding metrics. 
\end{thm}
\textbf{Proof.} We have only to prove $\tl H=\tl G^{\tl\sigma}$, which follows from \cite{an} (lemma 2.7). All the rest is an adaptation of \cite{Hel} (Ch. VII, thm 9.1) using what precedes. This completes the proof.\hfill$\blacksquare$\\[1.5mm]
\textbf{Acknowledgements} The author wishes to thank Francis Burstall for his suggestions which allowed him to improve the first version of this paper.

\vspace{0.2cm}
\noindent \textbf{Idrisse Khemar\\
TU Munich, Zentrum Mathematik (M8), Boltzmannstr.\,3, 85747 Garching, Germany\\
e-mail: khemar@math.jussieu.fr, khemar@ma.tum.de}


\begin{thebibliography}{111}
\bibitem{an} J. An, Z. Wang, \emph{On the realization of Riemannian symmetric spaces 
in Lie groups II} preprint arXiv: math/0504120.
\bibitem{besse} A.L. Besse, \emph{Einstein Manifolds},
Springer-Verlag, Berlin, Heidelberg, New York, 1987.
\bibitem{12} F.E. Burstall, F. Pedit, \emph{Harmonic maps via
Adler-Kostant-Symes Theory}, Harmonic maps and integrable systems,
A.P. Fordy, J.C. Wood (Eds.), Vieweg (1994), 221-272.
\bibitem{8} F.E. Burstall and J.H. Rawnsley, \emph{Twistor theory
for Riemannian Symmetric Spaces with applications to harmonic maps
 of Riemann Surfaces} Lect. Notes in Math., vol. 1424, Springer, 1990.
\bibitem{dieudonne} J. Dieudonné, \emph{\'Eléments d'analyse, Tome
2}, Gauthiers-Villars.
\bibitem{DPW} J. Dorfmeister, F. Pedit and H.-Y. Wu,
\emph{Weierstrass type representation of harmonic maps into
symmetric spaces}, Comm. in Analysis and Geometry, 6(4) (1998), p.
633-668.
\bibitem{HR1} F. Hélein and P. Romon, \emph{Hamiltonian stationary Lagrangian
surfaces in $\C^2$}, Comm. in Analysis and Geometry Vol. 10, N. 1, 2002, p.
79-126.
\bibitem{HR2} F. Hélein and P. Romon, \emph{Weierstrass representation
of Lagrangian surfaces in four dimensional spaces using spinors and
quaternions}, Comment. Math. Helv., 75 (2000), p. 668-680.
\bibitem{HR3} F. Hélein and P. Romon, \emph{Hamiltonian stationary
Lagrangian surfaces in Hermitian symmetric spaces}, in
\emph{Differential Geometry and Integrable Systems}, M. Guest,
R. Miyaoka, and Y. Ohnita, Editors-AMS, 2002.
\bibitem{Hel} S. Helgason, \emph{Differential geometry, Lie group
and symmetric spaces}, Academic Press, Inc., 1978.
\bibitem{jimenez} J.A. Jimenez, \emph{Riemannian 4-symmetric spaces},
Transactions of the American Mathematical Society, Vol.306, No.2. (Apr.,1988),
pp. 715-734. 
\bibitem{ki1} I. Khemar, \emph{Surfaces isotropes de $\oct$ et
syst\`{e}mes int\'{e}grables.},  J. Differential Geometry \textbf{79} (2008), no. 3, 479-516.
\bibitem{ki2} I. Khemar, \emph{Supersymmetric Harmonic Maps into
Symmetric Spaces}, Journal of Geometry and Physics 57 (2007) 1601-1630.
\bibitem{PS} A. Pressley and G. Segal, \emph{Loop Groups}, Oxford
Mathematical Monographs, Clarendon Press, Oxford, 1986.
\bibitem{tern} C.L. Terng, \emph{Geometries and Symmetries of Soliton
Equations and Integrable Elliptic Equations}, preprint
arXiv:math.DG/0212372.
\bibitem{wolf} J.A. Wolf and  A. Gray, \emph{Homogeneous spaces defined by Lie group automorphisms}. I, II,
J. Differential Geom. \textbf{2} (1968), 77-159.
\end{thebibliography}
\end{document}